\documentclass[11pt	]{article}
\usepackage[applemac]{inputenc}
\usepackage[T1]{fontenc}
\usepackage{natbib}
\usepackage[hidelinks]{hyperref}
\usepackage{color}
\usepackage{graphicx}
\usepackage[font = footnotesize]{subcaption}
\usepackage[font = footnotesize]{caption}
\usepackage{lscape}

\usepackage{enumitem}

\usepackage{amsmath}
\usepackage{amsfonts}
\usepackage{amsthm}
\newtheorem{remark}{Remark}[]
\newtheorem{lemma}{Lemma}[]
\newtheorem{theorem}{Theorem}[]
\newtheorem{proposition}{Proposition}[]
\newtheorem{corollary}{Corollary}[]

\usepackage{framed}
\usepackage[paper = a4paper, left = 25mm, right = 25mm]{geometry}
\renewcommand{\P}{\mathbb{P}}
\renewcommand{\d}{\mathrm{d}}
\newcommand{\wto}{\stackrel{\mathrm{w}}{\to}}

\DeclareMathOperator{\eqd}{\stackrel{\mathrm{d}}{=}}

\DeclareMathOperator{\wlim}{\mathrm{wlim}}
\newlength{\figwidth} 
\newcommand{\E}{\mathbb{E}}

\DeclareMathOperator{\R}{\mathbb{R}}

\newcommand{\N}{\mathbb{N}}

\DeclareMathOperator{\Var}{\mathrm{Var}}

\DeclareMathOperator{\Unif}{\mathrm{Unif}}
\newcommand{\change}[1]{{#1}}
\renewcommand{\[}{\left[}
\renewcommand{\]}{\right]}
\renewcommand{\(}{\left(}
\renewcommand{\)}{\right)}

\renewcommand{\d}{\mathrm{d}}
\newcommand{\iidsim}{\stackrel{\mathrm{i.i.d.}}{\sim}}

\newcommand{\Pois}{\mathrm{Pois}}

\newcommand{\ochi}{\overline{\chi}}
\newcommand{\oxi}{\overline{\xi}}
\newcommand{\oeta}{\overline{\eta}}
\newcommand{\okappa}{\overline{\kappa}}

\title{Critical cluster cascades}
\author{Matthias Kirchner\footnote{\href{mailto:matthias.kirchner@phnmsbern.ch}{matthias.kirchner@phnmsbern.ch}}}
\date{\today\footnote{to appear in {\it
Advances of Applied Probability}
in June 2023}}
\begin{document}
%
%



\maketitle
\begin{abstract}
We consider a sequence of Poisson cluster point processes on $\R^d$: At step $n\in\N_0$ of the construction, the cluster centers have intensity $c/(n+1)$ for some $c>0$, and each cluster consists of the particles of a branching random walk up to generation $n$---generated by a point process with mean 1. We show that this `critical cluster cascade' converges weakly, and that either the limit point process equals the void process (extinction), or it has the same intensity $c$ as the critical cluster cascade (persistence). We obtain persistence, if and only if the Palm version of the outgrown critical branching random walk is locally a.s.\ finite. This result allows us to give numerous examples for persistent critical cluster cascades.
\end{abstract} 


%

\section{Introduction}
In \cite{kallenberg17}, Chapter 13.10, Kallenberg summarizes work on `critical cluster stability' studied, e.g., in \cite{debes70}, \cite{ kallenberg77}, \cite{matthes78}, or \cite{liemant88}. In this context, one analyzes the effect of substituting all particles of some stationary point process on $\R^d$ by a `critical cluster', i.e., by shifted i.i.d. versions of a finite point process with mean total number of points equal to one. In particular, one identifies critical cluster distributions that allow nontrivial weak limits when such substitutions are iterated. If these iterations yield a nontrivial limit, the critical cluster (field) is called `stable'. Furthermore, one can show that the limit point process is invariant to further clustering, that is, it follows an equilibrium distribution. We closely follow the notation and line of argumentation as presented in \cite{kallenberg17} to study a similar limit construction:

We start with a homogeneous Poisson point process on $\R^d$. Its points form the roots of independent critical branching random walks. At each step of the construction and for each branching random walk, we either attach another generation of particles, or we delete it. This construction yields a sequence of Poisson cluster processes, where the cluster centers become fewer and fewer, and at the same time, the remaining clusters form `cumulative' critical branching random walks of more and more generations. The deleting and growing are balanced in such a way that the intensity stays unaffected. We call this sequence of point processes a `critical cluster cascade'. Note that in contrast to the iterated clustering construction in \cite{kallenberg17}, in our case a critical cluster cascade considers  more and more generations of the constructed critical branching random walks rather than `cousins' of higher and higher degrees.

For an overview of branching  random walks, see, e.g., \cite{shi16}. Due to almost sure extinction in the critical case, critical branching random walks are studied less frequently than supercritical branching random walks. In order to obtain interesting limiting properties, research on single critical branching random walks mostly conditions on survival. For instance, \cite{kesten95} studies maximal displacements of critical branching random walks on the real line with symmetric i.i.d. displacement and conditioned on survival. Critical branching random walks on integer lattices (of general dimension) with i.i.d. nearest neighbor displacement conditional on survival are studied in \cite{lalley11}; in the same setup, \cite{pekoez20} identifies second-order properties of the limiting distribution.

In contrast, we obtain nontrivial limits without conditioning on survival because critical cluster cascades involve the superposition of infinitely many critical branching random walks---rooted all over space. We show that if the generating critical cluster is sufficiently spread out, then persistence of the critical cluster cascade is possible:  in a way to be made precise, displacement outweighs extinction. We will show that persistence implies local integrability of the critical cluster cascade so that the limit process has the same intensity as each process in the critical cluster cascade sequence.  Furthermore, \change{in the persistent case, the limit process has an interesting structure: all particles are in a way explained by other particles and hence the process enjoys a kind of self-balancing structure. All explanatory power comes from within the system.}

\change{
The next section fixes terminology and introduces critical cluster cascades formally. In Section~\ref{sec:weak_convergence_and_persistence}, we prove weak convergence of critical cluster cascades and give a first persistence criterion. Section~\ref{sec:infinite_palm_tree} presents the infinite Palm tree, the limit of the Palm version of a cumulative critical branching random walk. Section~\ref{sec:persistence_theorem} gives the main result of the paper: the persistence theorem describes persistence of a critical cluster cascade in terms of local finiteness of the corresponding infinite Palm tree. We derive some corollaries that provide simple sufficient and necessary persistence conditions based on random-walk concepts. In Section~\ref{sec:examples}, we discuss examples. In Section~\ref{sec:outlook}, we formulate some open questions.  Appendix~\ref{appendix:proofs} collects the more technical proofs. In Appendix~\ref{appendix:figures}, we provide three figures that illustrate the various constructions.
}

\section{Model}\label{sec:model}
\change{Fix a dimension $d\in\N$. Let $\mathcal{N}_d$ be the set of locally finite counting measures on $(\R^d,\mathcal{B}(\R^d))$
equipped with the $\sigma$-algebra $\mathcal{A}_d$ generated by the sets $\{\mu\in\mathcal{N}_d: \mu B= k\}$ for $k\in\N_0$ and $B\in\mathcal{B}_b(\R^d)$. For any $\mu\in\mathcal{N}_d$, denote $\|\mu\|:=\mu\R^d$ 
and $\theta_x \mu := \mu(\cdot -x),\, x\in\R^d$. We call a random variable $\xi: (\Omega, \mathcal{F})\to(\mathcal{N}_d,\mathcal{A}_d) $ a \emph{point process}.
%
Any point process $\xi$ has an (possibly infinite) \emph{intensity measure} $\E \xi$ on $(\R^d,\mathcal{B}(\R^d))$.
A sequence of point processes $(\xi_n)$ \emph{converges weakly} to a point process $\xi_\infty$ if $\lim_{n\to\infty} \int f(x)\xi_n(\d x) = \int f(x)\xi_\infty(\d x)$ for all nonnegative continuous functions $f:\R^d\to\R_{\geq 0}$ with compact support. We write $\wlim_{n\to\infty} \xi_n:=\xi_\infty$. We also write $\wlim_{n\to\infty} X_n$ for the distributional limit of a sequence of univariate random variables $(X_n)\subset\R$). We call a sequence of point processes $(\xi_n)$ \emph{locally uniformly integrable} if $(\xi_n B)$ is uniformly integrable for all bounded Borel sets $B\subset\R^d$. For any sequence of random variables $(S_n)\subset\R^d$, we call $\sum \delta_{S_n}$ its \emph{occupation measure} and $\E\sum \delta_{S_n}$ its \emph{expected occupation measure}. We write $\mathcal{L}(X)$ for the distribution of a random variable $X$. For two measures on $\R^d$, $\rho_1$ and $\rho_2$, we define their \emph{convolution} as $\rho_1*\rho_2(\cdot):=\int\rho_1(\cdot -x)\rho_2(\d x)$. We denote the open ball with radius $r>0$ and center $x\in \R^d$ by $B_x^r$. Finally, $\lambda^d$ denotes $d$-dimensional Lebesgue measure.}\\

Let
$
\chi
$
be a point process with $\E\|\chi\| =1$ and $\Var\|\chi\| <\infty$. Furthermore, throughout the paper, we assume that $\chi$ is \change{simple (i.e., \change{$\P\{\chi\{x\}\leq 1,\ x\in\R^d\}=1$}) and diffuse (i.e., $\P\{\chi\{x\} = 0\} = 1$, $x\in\R^d$)}. We call such a $\chi$ a \emph{critical cluster}, and $(\chi^x)_{x\in\R^d}$, with $\chi^x:\eqd \theta_x\chi$ independent over $x\in\R^d$, a \emph{critical cluster field}. For $x\in\R^d$, set
\begin{equation}
\chi^x_0 := \delta_x,\quad \chi^x_k := \int\chi^y\chi^x_{k-1}(\d y),\quad k\in\N.\label{eq:branching_random_walk}
\end{equation}
Note that, by the diffuseness and simpleness assumption on the critical cluster, for all $x\in\R^d$, $(\|\chi^x_k\|)_{k}$ forms a critical Galton--Watson process with offspring distribution $\mathcal{L}(\|\chi\|)$, and  $(\chi^x_k)_{k}$ forms a \emph{critical branching random walk} generated by $\mathcal{L}(\chi)$ and rooted in $x$. \change{In particular, we denote the branching random walk rooted in zero, $(\chi^0_k)_k$, by $(\chi_k)_k$.}

Next, for $x\in\R^d$, define the \emph{cumulative branching random walk} $(\ochi^x_n)_{n}$ by
\begin{equation}
\overline{\chi}^x_n :=
\sum_{k = 0}^n\chi^x_k  ,\quad n\in\N_0,
\end{equation}
\change{so that $\ochi_n$ consists of the first $n+1$ generations of the branching walk $(\chi^x_k)_k$. We call the point $x$ (or the point processes $\ochi_0^x = \chi_0^x = \delta_x$) the \emph{root} of $\ochi^x_n$, and (the points measured by) $\chi^x_n$ the \emph{leafs} of $\ochi^x_n$.
As before, we denote $(\ochi_n)_n:=(\ochi^0_n)_n$.} Furthermore, let $\mu$ a Poisson process with constant intensity $c>0$, and, for $n\in\N_0$, \change{define its thinnings}
 \begin{equation}
 \mu_n(\d x):= 1\left\{U_x \leq \frac{1}{n+1}\right\}\mu(\d x),\quad \text{with } U_x\iidsim \Unif(0,1),\, x\in\R^d,\label{eq:thinnings}
 \end{equation}
 so that
 $(\mu_n)$ forms an a.s.\ nonincreasing sequence of Poisson processes on $\R^d$ with intensity sequence $(c/(n+1))_n$. \change{We call the particles measured by  $(\mu_n)$ \emph{immigrant points}.}

\begin{figure}[h]
\centering
\begin{subfigure}{\figwidth}
\includegraphics[width = \figwidth]{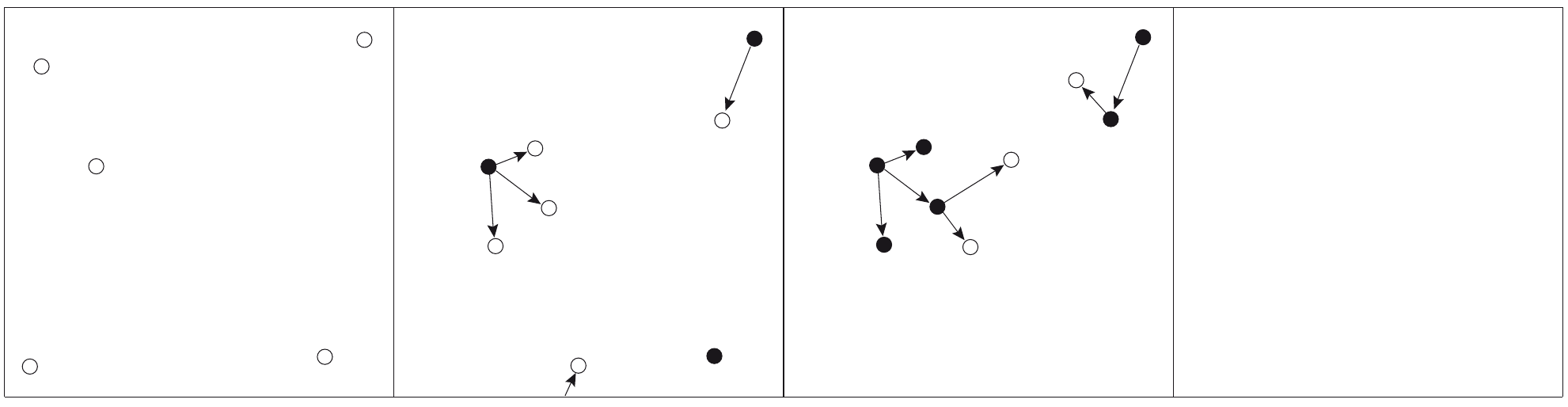}
\caption{Extinction: The distribution of the displacement of the cluster points is so concentrated around its respective origin that the points of each branching random walk tend to clump. At the same time, these clumps become sparser and sparser in space so that in the limit we ultimately obtain the void process.  Corollary~\ref{cor:nec} is helpful to find examples for critical cluster cascades that extinguish; see Section~\ref{example:d=1,2}.
}
\end{subfigure}
\begin{subfigure}{\figwidth}
\medskip
\includegraphics[width = \figwidth]{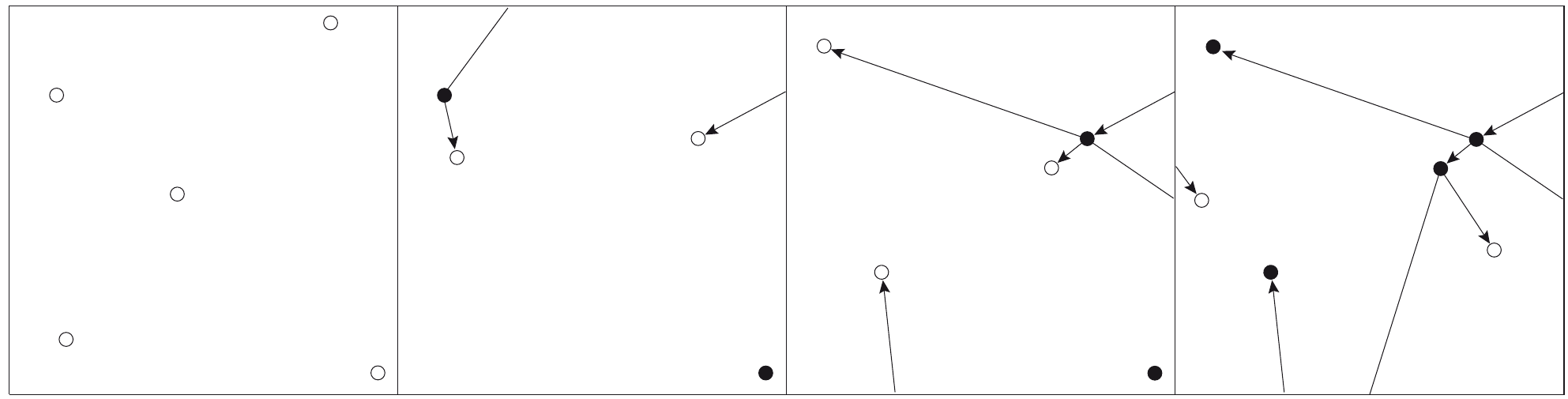}
\caption{Persistence: The distribution of the displacement of the clusters is so spread out that clumping of the branching random walks is avoided, and it will always remain possible to observe some points. Though, in the limit, in any finite set, we will not observe any immigrant points, anymore; see Remark~\ref{rmk:properties}\ref{rmk:properties:cluster_centers}. Theorem~\ref{thm:persistence} shows, that the intensity of the limit process is the same as in all of the component processes. Corollaries~\ref{cor:suff1} and~\ref{cor:suff2} are helpful to find examples for critical cluster cascades that persist; see Sections~\ref{ex:alphastable} and~\ref{example:dgeq5}.}
\end{subfigure}

\caption{Illustration of two (hypothetical) realizations of the first four components $(\oxi_0, \oxi_1, \oxi_2, \oxi_3)$ of two (different) critical cluster cascades $(\oxi_n)$ in $\R^2$ as given in \eqref{eq:critical_cluster_cascade}. In both cases, we start with cluster centers (or `immigrant points') $\oxi_0=\mu$ as given in \eqref{eq:thinnings}. At step $n$, we either thin an immigrant point $x$ together with all of its previous offspring $\ochi^x_{n-1}$, or we attach a further generation of i.i.d. clusters $\chi^y$ to each of the leaf points $y$ measured by $\chi^x_{n-1}$ (empty points in illustration of $\oxi_{n-1}$)---thus creating new leafs (empty points in illustration of $\oxi_{n}$). All clusters $\chi^y$ have an expected number of points equal to one.
That is, the n-th component $\oxi_n$ of a critical cluster cascade consists of the particles of branching random walks up to generation $n$ attached to the remaining immigrants $\mu_n$.
Theorem~\ref{thm:existence} shows that critical cluster cascades converge weakly to some limit point process. Lemma~\ref{lemma:monotonicity} shows that we expect to observe fewer and fewer branching random walks (in average) in any finite set. Theorem~\ref{thm:kappa} and Theorem~\ref{thm:persistence} give criteria for whether the limit is the a.s.\ void point process (`extinction') or not (`persistence').}
\label{fig:1}
\end{figure}

Finally, consider a sequence of cluster processes $(\oxi_n)$, where $\oxi_n$ has the immigrants $\mu_n$ as cluster-center process and each cluster consists of the particles of the first $n+1$ generations of a branching random walk. I.e., 
\begin{equation}
\overline{\xi}_n := \int\overline{\chi}_n^x\mu_n(\d x),\quad n\in\N_0.\label{eq:critical_cluster_cascade}
\end{equation}
We call $(\oxi_n)$ a  \emph{critical cluster cascade}. \change{Figure~\ref{fig:1} gives an illustration.}

\begin{remark}\label{rmk:properties}
The following properties of a critical cluster cascade are easy to establish:
\begin{enumerate}[label = \rm(\roman*)]
\item The immigrant points die out, i.e., $\lim_{n\to\infty}\mu_n B = 0$ a.s.\ for any bounded Borel set $B$, and $\wlim_{n\to\infty} \mu_n  $ is the void point process. \label{rmk:properties:cluster_centers}
\item The cumulative critical branching random walks $(\ochi_n^x)_n$ converge a.s.\ to totally finite point processes $\ochi_\infty^x$ whenever $\Var\|\chi\| >0$. And, in any case, $\E\|\ochi_n^x\| = n+1$.\label{rmk:properties:finiteness}
\item The critical cluster cascade has constant intensity
$\E\oxi_n  = {c}\E\|\ochi_n\|/(n+1)  \lambda^d \equiv c \lambda^d,\ n\in\N_0$. (Recall that $\lambda^d$ denotes $d$-dimensional Lebesgue-measure.)
\end{enumerate}
\end{remark}

\section{Weak convergence and persistence of the critical cluster cascade}
\label{sec:weak_convergence_and_persistence}
We will show that $(\oxi_n)$ converges to a weak limit process $\oxi_\infty$. If $\oxi_\infty$ is nontrivial, we say that \emph{the critical cluster cascade persists}. Otherwise, \change{ i.e., if $\P\{\|\oxi_\infty\| = 0\} = 1$}, we say that \emph{the critical cluster cascade extinguishes}. \change{As it turns out in Section~\ref{sec:persistence_theorem}, if the critical cluster cascade persists, then the limit process has the same intensity $c\lambda^d$ as each process of the critical cluster cascade.}

For $r>0$ and $n\in\N_0$, consider
\begin{equation}
\okappa_n^r:=\int 1\{\ochi^x_nB_0^r > 0\}\mu_n(\d x),
\end{equation}
i.e., 
 the number of cumulative critical branching random walks in $\oxi_n$ hitting $B_0^r$.
Obviously, $\okappa_n^r$ is the total mass of the Poisson process $\mu_n$ after independent thinnings. Therefore, $\okappa_n^r$ is Poisson distributed and we immediately obtain the following result.
\begin{lemma}
The sequence $(\okappa_n^r)_n$ is uniformly integrable.\label{lemma:ui}
\end{lemma}
\begin{proof}
Uniform integrability follows from $\Var  \okappa_n^r=  \E \okappa_n^r \leq \E \oxi_nB_0^r = c\lambda^d B_0^r < \infty,\ n\in\N_0$.
\end{proof}
Furthermore, the expectations of the sequence converge:
\begin{lemma}\label{lemma:monotonicity}
For all $r>0$, $\E\okappa_n^r<\infty,\ n\in\N_0,$ and the sequence
$(\E\okappa_n^r)_n$ is nonincreasing. Furthermore, the limit $p_{r}:= \lim_{n\to\infty}\E\okappa_n^r$ is strictly positive for all $r > 0$ if and only if the limit $p_{r_0}$ is strictly positive for some $r_0>0$.
\end{lemma}
\begin{proof}
See \hyperref[proof:lemma:monotonicity]{Appendix A}.
\end{proof}
Hence,
$\okappa_n^r\wto \Pois(p_r),\quad n\to \infty,\quad\text{for some $p_r\in[0,\infty)$}, 
$ and $
\lim_{n\to\infty}\P\{\oxi_nB_0^r =0\} = \lim_{n\to\infty}\P\{\okappa_n^r = 0\}=\exp\{-p_r\}.
$
Note that the arguments above (and in the proof of Lemma~\ref{lemma:monotonicity}) can be repeated to find the existence of $\lim_{n\to\infty}\P\{\oxi_nB =0\}$ for arbitrary bounded Borel sets $B\subset\R^d$. We conclude that the void probabilities of $(\oxi_n)$ converge. This suffices for the following result.
\begin{theorem}[Existence of weak limit]\label{thm:existence}
The sequence of point processes $(\oxi_n)$ converges weakly to a point process $\oxi_\infty$. 
\end{theorem}
\begin{proof}
Weak convergence follows from the convergence of the void probabilities; see, e.g., Theorem 2.2 in \cite{kallenberg17}. 
\end{proof}
\change{The next theorem gives a first criterion for persistence of a critical cluster cascade.}
\begin{theorem}[Persistence and extinction]\label{thm:kappa}
If $p_{r_0}=0$ for some $r_0>0$, then $(\oxi_n)$ extinguishes. If $p_{r_0}>0$ for some $r_0>0$, then $(\oxi_n)$ persists.
\end{theorem}
\begin{proof}
Clearly, $\P\{ \oxi_\infty B_0^r =0\} = \P\{\wlim_{n\to\infty}{\okappa_n^r} = 0\} = \exp(-p_r)$. Consequently, if $p_{r_0} = 0$ for some $r_0>0$ (and then, by Lemma~\ref{lemma:monotonicity}, for all $r>0$), then $\oxi_\infty B_0^r = 0 $ a.s.\ for all $r>0$. Similarly, $\P\{ \oxi_\infty B_0^r >0\} >0$ for $r>0$ if $p_{r_0} >0$ for some $r_0>0$.
\end{proof}

\change{Unfortunately, the persistence criterion in Theorem~\ref{thm:kappa} is only useful in very special cases. For instance, if we consider clusters $\chi$ without displacement; see Section~\ref{sec:clusters_without_displacements}. In order to find more convenient criteria for persistence, we study a Palm version of a single cumulative branching random walk $(\ochi_n)$:
}
\section{Infinite Palm tree}\label{sec:infinite_palm_tree}
\begin{figure}
\centering
\includegraphics[width = \figwidth]{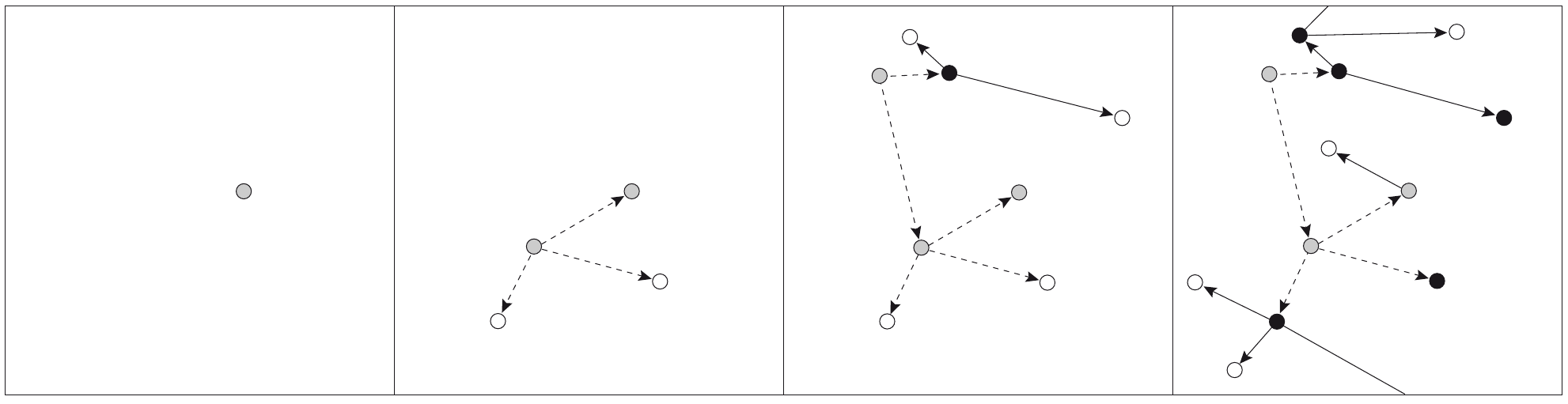}
\caption{Illustration of the first four steps $(\oeta_0^{(L_0)},\oeta_1^{(L_1)},\oeta_2^{(L_2)},\oeta_3^{(L_3)} )$ of a hypothetical realization of the forward/backward construction of an infinite Palm tree $\oeta_\infty$ in $\R^2$ in (the proof of) Lemma~\ref{lemma:fwbw}; see the recursion in~\eqref{eq:fwbw_recursion}. We start with a single point in zero (single grey point in first panel). The increments ($\in \{0,1\})$ of the Markov chain $(L_n)$ defined in \eqref{eq:L} determine whether to perform a (genealogical) forward step (first case in \eqref{eq:fwbw_recursion}) or a backward step (second case in \eqref{eq:fwbw_recursion}). In our illustration, we first realize two backward steps and then a forward step. In the first backward step,  we attach a possible parent point (in grey) to the zero point together with its sibling points (empty points); that is, the points connected by the dashed lines are a realization of the parent/siblings process $(\beta^0_{0,0}, \beta^0_{1,0})$ defined in \eqref{eq:parentsibling_shifted}. Together with the zero point, these four points constitute $\oeta_1^{(L_1)}$. We proceed with another backward step. This time, together with the parent of the earlier parent and its sibling (black point), we attach another generation of clusters $\chi^x$ to the siblingt $x$ so that the tree $\oeta^{(L_2)}_2$ consists of three generations. We refer to the foremost generation of points (empty points), $\eta^{(L_2)}_2$, as leaf points; see \eqref{eq:uetal}. In the following forward step, we attach a cluster $\chi^y$ to each `leaf point' $y$  (and thus generating new leafs and a new generation of points). Proposition~\ref{prop:direct_construction_of_infinite_palm_tree} gives a more direct construction of the limit object, the infinite Palm tree. This more direct construction is illustrated in Figure~\ref{fig:3}.}
\label{fig:2}
\end{figure}

The limit behavior of $(\overline{\xi}_n)$ is intimately related to the limit behavior of a Palm version $(\oeta_n)$ of the cumulative critical branching random walk $(\ochi_n)$. \change{The marginal distributions of $(\oeta_n)$ are determined by
\begin{equation}
\E f(\oeta_n) := \frac{1}{n+1} \E \int f\big(\ochi_n(\cdot + x) \big) \ochi_n(\d x),\quad \forall f:\mathcal{N}_d\to\R_0^+, \text{ measurable},\ n\in\N_0.\label{eq:oeta}
\end{equation}
The point process $\oeta_n$ corresponds to a version of $\ochi_n$ that is shifted in space in such a way that one of its points lies in zero. That is, just like $\ochi_n$,  $\oeta_n$ consists of $n+1$ generations, it has a root point (at some random location) and points of a `last' (possibly void) generation $n$ that we call \emph{leafs}; this genealogical structure of $\oeta_n$ will be made precise in the proof of Lemma \ref{lemma:fwbw} below.}
\change{
\begin{remark}
The following properties of $\oeta_n$ are easy to establish:
\begin{enumerate}[label = \rm(\roman*)]
\item There is a.s.\ exactly one point in zero: $\P\{\oeta_n\{0\} = 1\} =  \E \int 1\{\ochi_n\{x\} = 1 \} \ochi_n(\d x)/(n+1) = \E\|\ochi_n\|/(n+1) = 1.$ And, in particular, $\|\oeta_n\|>0$ a.s.
\item The distribution of the total number of points of the process $\oeta_n$ is $\P\{\|\oeta_n\| = k\} =\E1\{\|\ochi_n\|=k\}\|\ochi\|/(n+1)  = k\P\{\|\ochi_n\|=k\}/(n+1)$, $k\in\N_0$. I.e., the distribution of $\|\oeta_n\|$ is a size-biased version of the distribution of the total number of points of the original process $\ochi_n$.
\item The expected total number of points of $\oeta_n$ is $\E\|\oeta_n\| = \E \|\ochi_n\|^2/(n+1)$.
\end{enumerate}
\end{remark}

Lemma~\ref{lemma:fwbw} below shows that it is possible to give an a.s.\ nondecreasing construction of a sequence of point processes $(\oeta_n)$ such that its marginal distributions are determined by~\eqref{eq:oeta}. Consequently, it makes sense to define the random measure $\oeta_\infty:=\lim\oeta_n$. We call this limit measure the \emph{infinite Palm tree}. As it turns out in Proposition~\ref{prop:01}, finiteness of the {infinite Palm tree} in a neighborhood around zero is a 0--1 event. In the next section, Theorem~\ref{thm:persistence} will show that local finiteness of the infinite Palm tree is equivalent to persistence of the corresponding critical cluster cascade.

The construction of the infinite Palm tree in (the proof of) Lemma~\ref{lemma:fwbw} is similar to the `method of backward trees' in \cite{kallenberg77} and \cite{liemant81}. However, in contrast to this earlier work, we are not only interested in a Palm version of the particles of generation $n$ (i.e., of $\chi_n$) but of all generations \emph{up to} generation $n$ (i.e., of $\ochi_n=\sum_{k=0}^n\chi_k$).
%
%

In this construction and also later in the paper, we will make use of the Palm version of the tuple 
$(\chi_0,\chi_1)(\eqd(\delta_0, \chi))$, i.e., of generation 0 and generation 1 of our critical branching random walk $(\chi_n)$ defined in   \eqref{eq:branching_random_walk} and rooted in zero:
 \begin{equation*}
 \E f(\eta^{(1)}_{0}, \eta^{(1)}_{1}) := \E\int f\big(\chi_0(\cdot + x), \chi_1(\cdot + x)\big) \chi_1(\d x), \quad \forall f:\mathcal{N}^2_d\to\R_0^+, \text{ measurable}.
\end{equation*}
Note that $\eta^{(1)}_{1}$ is a shifted version of the cluster $\chi_1 (\eqd \chi)$ given it has a point in zero. I.e., $\eta^{(1)}_{1}$ consists of the point zero and its potential \emph{siblings}. The process $\eta^{(1)}_{0}$ measures the center or \emph{parent} of the same shifted version of the cluster $\chi_1$ given it has a point in zero. I.e., $\eta^{(1)}_{0}$ consists of exactly one point, the parent. We further consider an i.i.d. sequence of the same process, disregarding the point in zero:
\begin{equation}
(\beta_{0,n}, \beta_{1,n}):\eqd (\eta_0^{(1)}, \eta_1^{(1)}- \delta_0),\quad \text{ independently over }n\in\N_0.\label{eq:parentsibling}
\end{equation}
We call $(\beta_{0,n}, \beta_{1,n})$ \emph{parent/siblings processes}.
Finally, for all $n\in\N_0$, we consider the shifted versions of the parent/siblings processes:
\begin{equation}
\big(\beta^x_{0,n}, \beta^x_{1,n}):=\big(\beta_{0,n}(\cdot -x), \beta_{1,n}(\cdot - x)\big),\quad x\in\R^d.\label{eq:parentsibling_shifted}
\end{equation}
 Note that $\beta^x_{0,n}$ consists of the cluster center or parent of a shifted version of the cluster $\chi_1$ given it has a point in $x$, and $\beta^x_{1,n}$ consists of the (potentially) remaining points of the cluster (disregarding the point in $x$). The parent/siblings process will play a roll in the construction of $(\oeta_n)$: 
 
 We start the recursive tree construction with 
$\oeta_0:=\delta_0$, i.e., with a single point at zero. 
At each step, we either perform a forward step or a backward step. The decision between forward an backward step will depend on the realization of a specific Markov chain.
In the forward step, we add particles by attaching another generation of clusters $\chi^x$ to leaf points of the previous tree. This forward step corresponds to adding another generation of clusters in our standard cumulative branching random walk, where $\ochi_n = \ochi_{n-1} + \int\chi^x\chi_{n-1}(\d x)$. Genealogically speaking, we grow the tree in a forward direction.
In the backward step, we attach a parent/siblings process to the root of the previous tree together with a specific number of offspring generations of these siblings. I.e., we grow the tree backwards. This construction is illustrated in Figure~\ref{fig:2} and will be made precise in the following proof.}

\begin{lemma}[Forward/backward construction of infinite Palm tree]\label{lemma:fwbw}
There exists an a.s.\ nondecreasing sequence of point processes $(\oeta_n)$ such that, for all $n\in\N_0$, the marginal distribution $\mathcal{L}(\oeta_{n})$ is  given by \eqref{eq:oeta}. 
\end{lemma}
\begin{proof}
Let $n\in\N_0$. Given a generation-wise vector representation $(\chi_0,\chi_1,\dots, \chi_n)$ of the cumulative branching walk $\ochi_n$, we define the distribution of its {Palm version} $(\eta^{(l)}_0, \eta^{(l)}_1, \dots,\eta^{(l)}_n)$ {with respect to the $l$-th generation} (for $l\in\N_0$) by
\begin{equation}
\E f(\eta^{(l)}_0, \eta_1^{(l)}, \dots,\eta_n^{(l)}) := \E\int f\(\theta_{-x}{\chi}_0,\theta_{-x}{\chi}_1, \dots, \theta_{-x}{\chi}_n \)\chi_l(\d x), \quad \forall f:\mathcal{N}^{n+1}_d\to\R_0^+, \text{ measurable}.
\label{eq:uetal}
\end{equation}
Set $\oeta^{(l)}_n:=\sum_{k=0}^n\eta^{(l)}_k$. Obviously, $
\E f(\oeta_n^{(l)}) = \E\int f(\theta_{-x}\ochi_n)\chi_l(\d x),\ \forall f:\mathcal{N}_d\to\R_0^+, \text{ measurable}.$
We call $\eta_0^{(l)}$ the \emph{root} and $\eta_n^{(l)}$ the \emph{leafs} of $\oeta^{(l)}_n$. Clearly, we retrieve the Palm version $\oeta_n$ of $\ochi_n$ (see \eqref{eq:oeta}) by 
\begin{equation*}
\E f(\oeta_n) = \frac{1}{n+1}\E\int f(\theta_{-x}\ochi_n)\ochi_n(\d x) =
\frac{1}{n+1}\sum_{l=0}^n\E f(\oeta_{n}^{(l)}).
\end{equation*}
In other words,
\begin{equation}
 \oeta_n\eqd \oeta_n^{(U_n)},\label{eq:mixture}
\text{ for $U_n\sim\Unif\{0,1,\dots,n\}$, independent of $\oeta_n^{(l)},\ l = 0,1, \dots, n$.} 
\end{equation}
Define recursively a random sequence $(L_n)$ by
$
L_0 := 0
$ and, for $n\in\N_0$ and $l\in\{0,1,2,\dots, n-1\}$,
\begin{equation}
\P\[L_{n+1} = l | L_{n} =l\] := \frac{n+1-l}{n+2}\quad\text{and}\quad \P\[L_{n+1} = l+1| L_{n} =l\] := \frac{l+1}{n+2}.\label{eq:L}
\end{equation}
By construction, $L_{n+1} - L_{n} \in\{0,1\}$ a.s., $n\in\N_0$. Furthermore, one can show by induction that  $L_n\sim\Unif\{0,1,\dots, n\},\ n\in\N_0$, so that 
$\mathcal{L}(\oeta_n^{(L_n)})=\mathcal{L}(\oeta_n)$.

\change{We aim to show that $(\oeta_n^{(L_n)})$ can be chosen nondecreasing.} Let $\oeta_0^{(0)} := \delta_0$.
Because $L_{n+1}-L_{n}\in\{0,1\}$, it suffices to find pathwise nondecreasing construction steps from
$\oeta^{(l)}_n$ to $\oeta^{(l)}_{n+1}$ ({forward step}) and from $\oeta^{(l)}_n$ to $\oeta^{(l+1)}_{n+1}$ ({backward step}) for $l\in\{0,1,\dots, n\}$. Given the sequence $(L_n)$, we can then construct a nondecreasing sequence of random measures $(\oeta^{(L_n)}_n)$ with $\oeta^{(L_n)}_n\eqd \oeta_n$ (see \eqref{eq:mixture}), 
where in step $n$ we either make a forward or a backward step---depending on the realized step size of $L_n$.

\change{
So what is left to prove, is that we can construct the forward step and the backward step in such a way that in both cases we only add new points and do not remove old points. It might be helpful if the reader follows the proof together with Figure~\ref{fig:2}.}
\medskip\\
{\bf\nopagebreak Forward step: $\oeta_{n}^{(l)}\mapsto \oeta_{n+1}^{(l)}$}\\
We attach clusters $\chi^y$ to the leafs $\eta_{n}^{(l)}$ of $\oeta_{n}^{(l)}$ and find for all $f:\mathcal{N}_d\to\R_0^+$
\begin{align*}
\E f\(\oeta_n^{(l)} + \int \chi^y \eta_n^{(l)}(\d y)\)
&\overset{\eqref{eq:uetal}}{=} \E\int f\left(\theta_{-x}\ochi_{n} +\theta_{-x}\int\chi^y\chi_n(\d y)\)\chi_l(\d x)\\
&= \E \int f\(\theta_{-x}\ochi_{n+1}\)\chi_{l}(\d x)\\
&= \E f\(\eta_{n+1}^{(l)}\),\quad l = 0, 1,\dots, n.
\end{align*}
So we have shown that the forward step yields the desired distribution.\\
\noindent {\bf Backward step: $\oeta_{n}^{(l)}\mapsto \oeta_{n+1}^{(l+1)}$}\\
Also the backward step is conceptually simple. To the root $\eta_{0}^{(l)}$ of $\oeta_{n}^{(l)}$, we attach a parent point and potential sibling points. Furthermore, to each of these siblings we attach the first $n$ generations of a branching random walk. However, the notation is quite involved. First of all, we 
remind the reader of the \emph{parent/siblings processes} $({\beta}^x_{0,n}, \beta^x_{1,n})$ (independent over $n\in\N_0$)
from \eqref{eq:parentsibling_shifted}.

Given $\oeta_n^{(l)}$, we attach the {parent}-process $\beta^y_{0,n}$ to the root point $y$ (measured by $\eta_{0}^{(l)}$) of $\oeta_n^{(l)}$ which gives the new root (process)
\begin{equation}
\int\beta^y_{0,n}\eta_{0}^{(l)}(\d y) \label{eq:root}.
\end{equation}
 Furthermore, to each of the siblings $\int\beta_{1,n}^z\eta_{0}^{(l)}(\d z) $ of the old root $\eta^{(l)}_{0}$, we attach shifted (and independent) versions of the cumulative branching random walk $\ochi_n$ 
 \begin{equation}
 \int \int\ochi^y_n\beta_{1,0}^z(\d y)\eta_{0}^{(l)}(\d z)\label{eq:siblings}.
 \end{equation}
 Summarizing, the backward step consists of attaching the new root from \eqref{eq:root} and the potential siblings together with their offspring from \eqref{eq:siblings} to $\oeta_n^{(l)}$. One can show that the result of this backward step does indeed have distribution $\mathcal{L}(\oeta_{n+1}^{(l+1)})$. In fact, for $f:\mathcal{N}_d\to\R_0^+, \text{ measurable}$,
\begin{align}
& \E f\(\oeta_n^{(l)}  + \int\beta^y_{0,n}\eta_{0}^{(l)}(\d y) + \int\int\ochi^z_n\beta^y_{1,n}(\d z)\eta_{0}^{(l)}(\d y)\)\nonumber\\
&= \E f\(\oeta_n^{(l)} + \int\(\beta^y_{0,n} + \int\ochi_n^z\beta^y_{1,n}(\d z)\)\eta_{0}^{(l)}(\d y)\)\nonumber\\
&\stackrel{\eqref{eq:uetal}}{=} \E \int f\(\theta_{-x}\[\ochi_n + \int\(\beta^y_{0,n} + \int \ochi_n^z\beta^y_{1,n}(\d z)\)\chi_0(\d y)\]\)\chi_l(\d x).\label{eq:lastline_alt}
\end{align}
Noting that $\chi_0=\delta_0$, and denoting $(\tilde{\eta}^{(1)}_0,\tilde{\eta}_1^{(1)})$ respectively $\tilde{\chi}_1$ as independent copies of $({\eta}^{(1)}_0,{\eta}_1^{(1)})$ respectively ${\chi}_1$, and identifying $\ochi^0_n$ with $\ochi_n$, we obtain
\begin{align}
\eqref{eq:lastline_alt}&= \E \int f\(\theta_{-x}\[\ochi_n +\beta^0_{0,n} +  \int \ochi_n^z\beta^0_{1,n}(\d z)\]\)\chi_l(\d x)\nonumber\\
&\stackrel{\eqref{eq:parentsibling}}{=} \E \int f\(\theta_{-x}\[\ochi_n +\tilde{\eta}^{(1)}_0 +  \int \ochi_n^z(\tilde{\eta}_1^{(1)}-\delta_0)(\d z)\]\)\chi_l(\d x)\nonumber\\
&=\E \int f\(\theta_{-x}\[\tilde{\eta}_{0}^{(1)} +  \int \ochi_n^z\tilde{\eta}_1^{(1)}(\d z)\]\)\chi_l(\d x)\nonumber\\
&\stackrel{\eqref{eq:uetal}}{=} \E \int \int f\(\theta_{-x}\[\theta_{-y}\tilde{\chi}_0 +  \int \ochi_n^z\theta_{-y}\tilde{\chi}_1(\d z)\]\)\chi_l(\d x)\tilde{\chi}_1(\d y)\nonumber\\
&= \E \int \int f\(\theta_{-(x+y)}\[\delta_0 +  \int \ochi_{n}^z \tilde{\chi}_1(\d z)\]\)\chi_l(\d x)\tilde{\chi}_1(\d y)\nonumber\\
&= \E \int \int f\(\theta_{-(x+y)} \ochi_{n+1} \)\chi_l(\d x)\tilde{\chi}_1(\d y)\nonumber\\
&= \E \int \int f\(\theta_{-x} \ochi_{n+1} \)\theta_{-y}\chi_l(\d x)\tilde{\chi}_1(\d y)\label{eq:lastline}.
\end{align}
We have 
$\int\theta_{-y}\chi_l(\d x)\tilde{\chi}_1(\d y) =  \chi_{l+1}(\d x).
$ So we may conclude the calculation by
\begin{align*}
\eqref{eq:lastline} = \E \int f\(\theta_{-x}(\ochi_{n+1})\)\chi_{l+1}(\d x) 
= \E  f\(\oeta_{n+1}^{(l+1)}\),
\end{align*}
and also the backward step gives the desired distribution.
We summarize: given the sequence $(L_n)$ as defined in \eqref{eq:L}, we set $\oeta_0 := \delta_0$, and  define recursively
\begin{align}
&\oeta_{n+1}^{(L_{n+1})}\nonumber\\
&:=
\begin{cases}
\oeta_n^{(L_{n})} + \int \chi^y \eta_{n}^{(L_{n})}(\d y), & \text{if }L_{n+1} = L_{n} \quad\text{(forward step),}\\
\oeta_n^{(L_{n})}  + \int\beta^y_{0,n}\eta_{0}^{(L_{n})}(\d y) + \int\int\ochi^z_n\beta^y_{1,n}(\d z)\eta_{0}^{(L_{n})}(\d y),&\text{if } L_{n+1} = L_{n} +1  \quad\text{(backward step)}.
\end{cases}\label{eq:fwbw_recursion}
\end{align}
Obviously, the sequence of measures $(\oeta_n^{(L_{n})})$ is a.s.\ nondecreasing. In addition, the calculations given above \eqref{eq:mixture} show that for all $n\in\N_0$ the law of $\oeta_n^{(L_{n})}$ coincides with the law of $\oeta_n$. So we may choose $(\oeta_n^{(L_{n})})$ as candidate for the a.s\ nondecreasing version of $(\oeta_n)$.
\end{proof}
Because of Lemma~\ref{lemma:fwbw}, it makes sense to define the a.s.\ limit random measure $\oeta_\infty := \lim_{n\to\infty}\oeta_n$. We call $\oeta_\infty$ the \emph{infinite Palm tree}. Unfortunately, the construction in the proof of Lemma~\ref{lemma:fwbw} is not suited for further analysis of the infinite Palm tree. The next proposition will provide a more suitable representation of  (the distribution of) $\oeta_\infty$.

\begin{figure}
\centering
\includegraphics[width = \figwidth]{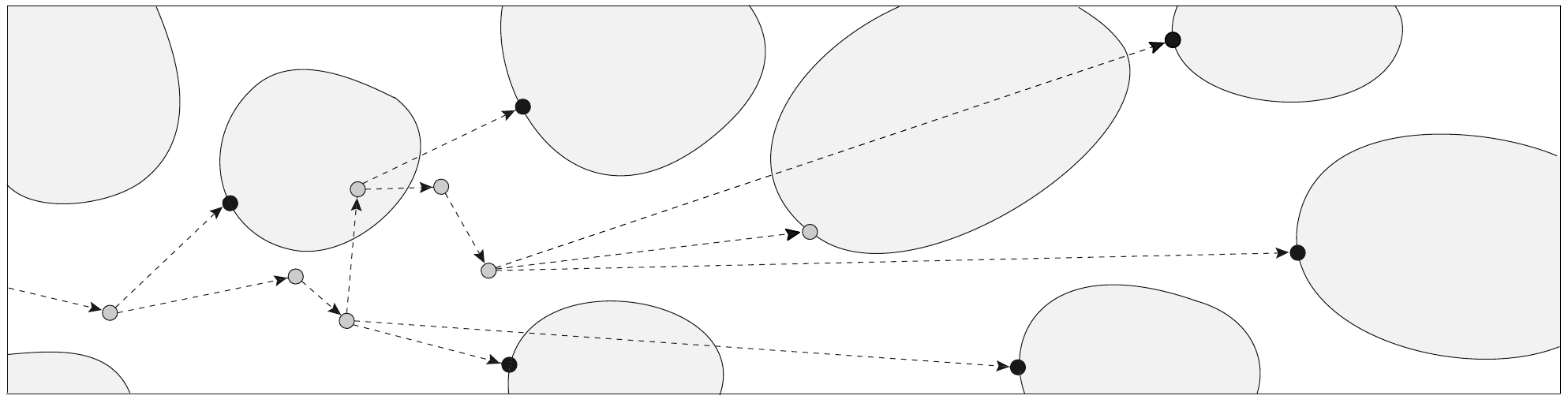}
\caption{Illustration of a hypothetical realization of the direct construction of the infinite Palm tree given in  Proposition~\ref{prop:direct_construction_of_infinite_palm_tree}. The grey points refer to the first seven values of the `infinite backward spine' random walk $(\zeta^-_n)$ of parents; see \eqref{eq:rw}. These points are generated by attaching a parent point $\zeta^-_1$ to the zero point $\zeta^-_0$ (grey point in center), a parent point $\zeta^-_2$ to the point $\zeta^-_1$, and so on. Note that these grey points correspond to the grey points in Figure~\ref{fig:2}. At each step, together with the parent point $\zeta^-_{n+1}$ of $\zeta^-_{n}$, we attach potential sibling points of $\zeta^-_{n}$ measured by $\beta^{\zeta^-_{n}}_{1,n}$ (black points). I.e., the dashed arrows indicate realizations of the parent/siblings processes $(\beta^x_{0,n}, \beta^x_{1,n})$; see \eqref{eq:parentsibling_shifted}. To each of these sibling points as well as to the zero point, we attach independent outgrown (however a.s.\ finite; see Remark~\ref{rmk:properties}\ref{rmk:properties:finiteness}) cumulative branching random walks $\ochi_{\infty, n}^x$ as given in \eqref{eq:outgrown} (shaded potato-like areas).}\label{fig:3}
\end{figure}

\change{We will prove that the construction described in the proof of Lemma~\ref{lemma:fwbw} and illustrated in Figure~\ref{fig:2} will involve infinitely many backward steps and infinitely many forward steps with probability 1. So, an alternative way to construct the infinite palm tree would be to construct an \emph{infinite backward spine} $(\zeta^-_n)_n$ of parents by setting
$\zeta^-_0:= 0$ and recursively attaching parent points (see \eqref{eq:parentsibling_shifted})
\begin{equation}
\quad \zeta^-_{n+1}:=\zeta^-_{n} + \int x\beta^{\zeta^-_{n}}_{0,n}(\d x), \quad n\in\N_0.\label{eq:rw}
\end{equation}
It is easy to show that $(\zeta^-_n)$ is a random walk with step-size distribution $\rho^-:=\E\chi(-\cdot)$; see \eqref{eq:rho-}. 

The infinitely many forward steps of the construction in the proof of Lemma~\ref{lemma:fwbw} and illustrated in Figure~\ref{fig:2} lead to \emph{outgrown cumulative branching random walks} defined by $\ochi_\infty := \lim_{n\to\infty} \ochi_n
$ and
\begin{equation}\label{eq:outgrown}
\ochi^x_{\infty,n}:\eqd \ochi_\infty(\cdot -x ) \quad \text{independently over $(n,x)\in\N_0\times \R^d$.}
\end{equation} 

Intuitively, we might want to attach these outgrown cumulative branching random walks directly to each of the points $(\zeta^-_n)_n$ in the infinite backward spine. However, the position of the parent $\zeta^-_{n+1}$ of the point $\zeta^-_{n}$ and the position (or number!) of its siblings are in general not independent. This is why, at the backward step, together with the parent point $\zeta^-_{n+1}$ (measured  by $\beta_{n,0}^{\zeta_n}$), we jointly have to model the sibling points (measured  by $\beta_{n,1}^{\zeta_n}$) of the point $\zeta^-_{n}$ with the \emph{parent/sibling processes} defined in \eqref{eq:parentsibling}. That is, we attach outgrown branching random walks	 to each of the potential siblings $\beta_{1,n}^{\zeta^-_{n}}$ of $\zeta^-_{n}$:
$$
\int\ochi^{x}_{\infty,n} B_0^r\beta_{1,n}^{\zeta^-_{n}}(\d x),\quad n\in\N.
$$
And finally, we attach the outgrown branching random walk $\ochi^0_{\infty,0}$ to the point zero (i.e.\ to $\zeta_0^-$).  
Proposition~\ref{prop:direct_construction_of_infinite_palm_tree} below summarizes this alternative construction of the infinite Palm tree; Figure~\ref{fig:3} gives an illustration.
}

\begin{proposition}[Direct construction of infinite Palm tree]\label{prop:direct_construction_of_infinite_palm_tree} With the notation from above, we have that the random measure $\oeta_\infty(\cdot)$ is equal in distribution to the random measure
 \begin{align}
 \ochi^0_{\infty,0}(\cdot)  + \sum_{n =0 }^\infty\(\int\ochi^{x}_{\infty,n} (\cdot) \beta_{1,n}^{\zeta^-_{n}}(\d x)  +  \beta^{\zeta^-_{n}}_{0,n}(\cdot) \).\label{eq:altpalmtree_measure}
 \end{align}
 \end{proposition}
\begin{proof}
See \hyperref[proof:prop:direct_construction_of_infinite_palm_tree]{Appendix A}.
\end{proof}
In Section~\ref{subsec:poisson}, we show that for $\chi$ a Poisson process, \eqref{eq:altpalmtree_measure} simplifies to attaching outgrown branching random walks directly to the backward infinite spine. \change{Next, we use the representation of the infinite Palm tree in Proposition~\ref{prop:direct_construction_of_infinite_palm_tree} to prove the following result.}

\begin{proposition}\label{prop:01}
The event $\{\oeta_\infty B_0^r < \infty\}$ has either probability 0 or probability 1.
\end{proposition}
\begin{proof}
From Proposition~\ref{prop:direct_construction_of_infinite_palm_tree}, it suffices to show that 
 \begin{align}
 \ochi^0_{\infty,0}B_0^r  + \sum_{n =0 }^\infty\(\int\ochi^{x}_{\infty,n} B_0^r\beta_{1,n}^{\zeta^-_{n}}(\d x)  +  \beta^{\zeta^-_{n}}_{0,n}B_0^r\)<\infty \label{eq:altpalmtree} 
 \end{align}
 is a 0--1 event. To that aim, note that all point processes involved in the random measure \eqref{eq:altpalmtree_measure} are a.s.\  finite. (We have $\|\ochi^{x}_{\infty,n}\|<\infty$ because $\Var\|\chi\| >0$ by assumption; see Remark~\ref{rmk:properties}\ref{rmk:properties:finiteness}). Consequently, \eqref{eq:altpalmtree} holds if and only if 
\change{
\begin{equation}
\int\ochi^{x}_{\infty,n} B_0^r\beta_{1,n}^{\zeta^-_{n}}(\d x)  +  \beta^{\zeta^-_{n}}_{0,n}B_0^r>0 \text{ for finitely many $n\in\N_0$.}
\label{eq:summands}
\end{equation}}This is a 0--1 event which follows from a general version of the Hewitt--Savage 0--1 law on exchangeable sequences as formulated, e.g., in Theorem 3.15 of \cite{kallenberg02}. Indeed, in our setting, the `infinite sequence of i.i.d.\ random elements' from the cited theorem is the sequence 
\change{
\begin{equation}
\Big\{ \((\beta_{0,n},  \beta_{1,n}), (\ochi^{x}_{\infty,n})_x\) \Big\}_{n\in\N_0}.\label{eq:seq}
\end{equation}
The random walk $(\zeta^-_n)$ as well as all measures used in event \eqref{eq:summands} and therefore the event itself can be determined from this sequence. Furthermore, for any $m\in\N$, the order of the first $m$ values of the random walk $(\zeta^-_n)$ defined in \eqref{eq:rw} does \emph{not} affect the summands in \eqref{eq:summands} with indices larger than $m$. Consequently, the event \eqref{eq:summands} does not depend on the order of the first $m$ elements of the sequence \eqref{eq:seq}. So, the event~\eqref{eq:summands} belongs to the exchangeable $\sigma$-field generated by the sequence \eqref{eq:seq}.}
\end{proof}

\change{
We have shown that the infinite Palm tree is either a.s.\  locally finite or a.s.\  locally infinite. Next, we will show that in fact, local finiteness of the infinite Palm tree is equivalent to persistence of the corresponding critical cluster cascade. The proof depends on truncation of measures.} For any measure $\mu$ on $\R^d$, we set
\begin{align}
\mu^{r,k}(\d x) :=1\{\mu B_{x}^r\leq k\} \mu(\d x),\quad r,k>0.\label{eq:truncation}
\end{align} 
The next lemma shows that the truncation $\ochi_n^{r,k}$ of the cumulative critical branching random walk $\ochi_n$ is closely related to its Palm version $\oeta_n$: 
\begin{lemma}\label{lemma:palmtruncation}
Let $(\ochi_n)$ a cumulative critical branching random walk, $(\oeta_n)$ the corresponding sequence of Palm versions from \eqref{eq:oeta}, and $\oeta_\infty$ its limit, the infinite Palm tree. Then the following results hold:
\begin{enumerate}[label = \rm(\roman*)]
\item For $n\in\N_0$, ${\E\|\ochi_n^{r,k}\|}/{(n+1)} = \P\{\oeta_nB_0^r\leq k\}$.\label{lemma:palmtruncation:truncation}
\item The sequence $({\E\|\ochi_n^{r,k}\|}/{(n+1)})_n$ is nonincreasing in $n$, and $\lim_{n\to\infty}{\E\|\ochi_n^{r,k}\|}/{(n+1)} =  \P\{\oeta_\infty B_0^r\leq k\}$.\label{lemma:palmtruncation:limit}
\item The double limit $\lim_{k\to\infty}\lim_{n\to\infty}{\E\|\ochi_n^{r,k}\|}/{(n+1)}$ exists and equals  $\P\{\oeta_\infty B_0^r<\infty\}$.\label{lemma:palmtruncation:doublelimit}
\end{enumerate}
\end{lemma}
\begin{proof}\label{proof:lemma:palmtruncation} \ref{lemma:palmtruncation:truncation}
 Let $(\oeta_n)$ be the increasing Palm tree from Lemma~\ref{lemma:fwbw} and recall the truncation notation from \eqref{eq:truncation}. Then
\begin{align*}
\frac{1}{n+1}\E\|\ochi_n^{r,k}\|
& = \frac{1}{n+1}\E \int 1\{\ochi_n B_x^r\leq k\}\ochi_n(dx)\\
& = \frac{1}{n+1}\E \int 1\{\theta_{-x}\ochi_n B_0^r\leq k\}\ochi_n(dx)\\
& = \E  1\{\oeta_n B_0^r\leq k\}\\
& = \P\{\oeta_n B_0^r\leq k\},\quad r,k>0.
\end{align*}
\ref{lemma:palmtruncation:limit} Monotonicity follows from \ref{lemma:palmtruncation:truncation} as  $(\P\{\oeta_n B_0^r\leq k\})_n$ is nonincreasing in $n$ (because, by Lemma~\ref{lemma:fwbw}, $ (\oeta_n B_0^r)_n$ can be chosen a.s\ nondecreasing). For the limit, note that
\begin{align}
\lim_{n\to\infty}\frac{\E\|\ochi_n^{r,k}\|}{n+1}=\lim_{n\to\infty}\P\{\oeta_n B_0^r\leq k\} = \P\{\oeta_\infty B_0^r\leq k\}.\label{eq:altpalmtreeeq}
\end{align}
\ref{lemma:palmtruncation:doublelimit} Since $\P\{\oeta_\infty B_0^r\leq k\}$ is nondecreasing in $k$, we may take the limit in with respect to $k$ in \eqref{eq:altpalmtreeeq} and find
\begin{equation*}
\lim_{k\to\infty}\lim_{n\to\infty}\frac{\E\|\ochi_n^{r,k}\|}{n+1} = \P\{\oeta_\infty B_0^r< \infty\}.
\end{equation*}
\end{proof}

\change{
Note that $\E \ochi_n^{r,k}/(n+1) = \E \|\ochi_n^{r,k}\|/\E \|\ochi_n\|\in[0,1]$ is a measure for `clumping' of the particles measured by $\ochi_n$. E.g., if the limit in
Lemma~\ref{lemma:palmtruncation}\ref{lemma:palmtruncation:limit} equals zero for all $k$ (and then also for the double limit in Lemma~\ref{lemma:palmtruncation}\ref{lemma:palmtruncation:doublelimit}), then $(\ochi_n)$ exhibits strong clumping as `most of the points' of $\ochi_\infty$ have more than $k$ points in their $r$-neighborhood and therefore become truncated. So Lemma~\ref{lemma:palmtruncation} connects the behavior of the infinite Palm tree around zero with clumping of the cumulative branching random walk $(\ochi_n)$.
}
\section{Criteria for persistence and extinction}\label{sec:persistence_theorem}
\change{We present the main theorem of the paper: Persistence of a critical cluster cascade is equivalent to local finiteness of the infinite Palm tree. Furthermore, the limit process of a persistent critical cluster cascade necessarily has the same intensity as the component processes.} 
\change{\begin{theorem}[Persistence of critical cluster cascades]\label{thm:persistence}
The following are equivalent:
\begin{enumerate}[label = (\roman*)]
\item[\rm(i)] The critical cluster cascade $(\oxi_n)$ persists, i.e., for all $r>0$, $\P\{\oxi_\infty B_0^r>0\} (= \lim_{n\to\infty}\P\{\okappa_n^r>0\})> 0$.\label{thm:persistence:persistence}
\item[\rm (i*)] For some $r_0>0$, $\P\{\oxi_\infty B_0^{r_0}>0\} (= \lim_{n\to\infty}\P\{\okappa_n^{r_0}>0\})> 0$.\label{thm:persistence:persistence*}
\item[\rm (ii)]The infinite Palm tree $\oeta_\infty$ of the outgrown branching random walk $\ochi_\infty$ is a.s.\  locally finite, i.e., $\P\{\oeta_\infty B_0^r< \infty\} = 1$ for all $r>0$. 
\label{thm:persistence:palmtree}
\item[\rm (ii*)]  For some $r_0>0$, $\P\{\oeta_\infty B_0^{r_0}< \infty\} = 1$. 
\label{thm:persistence:palmtree*}
\item[\rm (iii)] For all $r > 0$, $ \E\oxi_\infty B_0^r= \lim_{n\to\infty}\E\oxi_n B_0^r=c\lambda B_0^r$ for all $r>0$. I.e., the critical cluster cascade $(\oxi_n)$ is locally uniformly integrable, and
the limit process $\oxi_\infty$ has the same intensity $c(>0)$ as $\oxi_n$ for $n\in\N_0$. \label{thm:persistence:ui}
\item[\rm (iii*)] For some $r_0 >0$, $ \E\oxi_\infty B_0^{r_0}= \lim_{n\to\infty}\E\oxi_n B_0^{r_0}=c\lambda B_0^{r_0}$.\label{thm:persistence:ui*}
\end{enumerate}
\end{theorem}}

\begin{proof}[Proof of Theorem~\ref{thm:persistence}]
We first prove {(iii)}$\ \Rightarrow\ ${(i)}$\ \Rightarrow\ ${(ii)}$\  \Rightarrow\ ${(iii)}:

{(iii)} $\Rightarrow\ ${(i)} is trivial as $c>0$ by definition, and therefore, if {(iii)} holds, the limit process cannot be a.s.\ void. For {(i)} $\Rightarrow$ {(ii)}, we observe that
\begin{align*}
\P\{\oeta_n B_0^r\leq k\} 
&= \frac{1}{n+1}\E\|\ochi_n^{r,k}\|\lambda^dB_0^r \\
&= \E\int (\ochi_n^x)^{r,k}B_0^r \mu_n(\d x)\\
&\geq \E \oxi_n^{r,k}B_0^r\\
& =\E \int_{B_0^r}1\{\oxi_nB^r_x\leq k \}\oxi_n(\d x)\\
&\geq \E\int_{B_0^r}1\{\oxi_nB^{2r}_0\leq k \}  \oxi_n(\d x)\\
&=\E 1\{\oxi_nB^{2r}_0\leq k \}  \oxi_nB_0^r.
\end{align*}

Taking $\liminf_{n\to\infty}$ on both ends of the inequality and applying Fatou's lemma on the right hand side, we obtain $\P\{\oeta_\infty B_0^r\leq k\} \geq \E 1\{\oxi_\infty B^{2r}_0\leq k \} \oxi_\infty B_0^r$, and, after letting $k\to\infty$, 
\begin{equation*}
\P\{\oeta_\infty B_0^r< \infty\} \geq \E  \oxi_\infty B_0^r \geq \P\{\oxi_\infty B_0^r>0\}  >0.
\end{equation*}
So statement {(ii)} follows because $\{\oeta_\infty B_0^r< \infty\}$ is a 0--1 event by Proposition~\ref{prop:01}. For(ii) $\Rightarrow$(iii), consider
\begin{equation*}
\oxi_n^{(r,k)}:= \int(\ochi^x_{n})^{r,k}B_0^r \mu_{n}(\d x) \leq k\okappa_n^r.
\end{equation*}
 From  Lemma~\ref{lemma:ui}, we find that that the sequence $(\oxi_n^{(r,k)}B_0^r)_n$ is uniformly integrable  for all $r,k >0$. Therefore,
\begin{align*}
\E\oxi_\infty B_0^r &=\E\wlim_{{n}\to\infty}\int\ochi^x_{n}B_0^r\mu_{n}(\d x) \\
& \geq \E\wlim_{n\to\infty}\int(\ochi^x_{n})^{r,k}B_0^r\mu_{n}(\d x) \\
&=\lim_{n\to\infty}\E\int(\ochi^x_{n})^{r,k}B_0^r\mu_{n}(\d x) \\
&=  c\lambda^dB_0^r\lim_{n\to\infty}\frac{E\|\ochi^{r,k}_{n}\|}{n+1}\\
&\to  c\lambda^dB_0^r\P\{\oeta_\infty B_0^r < \infty\}.
\end{align*}
Here, we use uniform integrability in the third step and Lemma~\ref{lemma:palmtruncation} in the last step, where we let $k\to\infty$. 
So, we obtain
\begin{equation*}
 c\lambda^dB_0^r= \lim_{n\to\infty}\E\oxi_{n}B_0^r\leq\E\oxi_\infty B_0^r \leq \liminf_{n\to\infty}\E\oxi_{n} B_0^r =c\lambda^dB_0^r.
\end{equation*}
Thus, uniform integrability of $(\E\oxi_{n}B_0^r)_n$ follows from the convergence of its means; see, e.g., Lemma 4.11 in \cite{kallenberg02}. 

\change{Summarizing, we have now shown that (i)$ \Leftrightarrow$ (ii) $\Leftrightarrow$ (iii). Note that (i*) $\Leftrightarrow$ (ii*) $\Leftrightarrow$ (iii*) can be proven along exactly the same lines. Furthermore, (i*) $\Leftrightarrow$ (i) by the second statement in Lemma~\ref{lemma:monotonicity}.}
\end{proof}
From Theorem~\ref{thm:persistence}(ii), we obtain simple sufficient conditions for persistence of the critical cluster cascade based on the expected occupation measures $U$ and  $U^-$ of random walks generated by the probability distributions $\rho:=\E\chi$ and $\rho^-:=\E\chi(-\cdot)$, as well as on the characteristic function $\hat{\rho}$ of $\rho$:

  \begin{corollary}[Sufficient persistence condition]\label{cor:suff1}
  Let $\Var\|\chi\|\in(0,\infty)$ and let $U$ respectively $U^-$ be the expected occupation measure of a random walk with step-size distribution $\E\chi$ respectively $\E\chi(-\cdot)$. If the convolution $(U *U^-)B_0^r:= \int U B^r_{-x} U^-(\d x)<\infty$ for some $r>0$, then the critical cluster  cascade $(\oxi_n)$ generated by $\mathcal{L}(\chi)$ persists.  
  \end{corollary}  
\begin{proof}
  Lemma~\ref{lemma:bound} below shows that, for all $r>0$, $(U*U^-)B_0^r<\infty$ implies almost sure finiteness of $\oeta_\infty B_0^{r/2}$ so that $\P\{\oeta_\infty B_0^{r/2}<\infty\} =1$ which, by Theorem~\ref{thm:persistence}, is equivalent to  persistence of the critical cluster cascade.
 \end{proof}
\begin{lemma}\label{lemma:bound}
Let $\Var\|\chi\|\in(0,\infty)$. Then  $\E\oeta_\infty B_0^r<\infty$ if $(U*U^-)B_0^{2r}<\infty$.
\end{lemma}

\begin{proof}
See \hyperref[proof:lemma:bound]{Appendix A}.
\end{proof}

Note that the case $\Var\|\chi\|=0$ will be treated as a special example in Section~\ref{subsec:1ptclusters}.
Also note that the random-walk based necessary condition in Corollary~\ref{cor:suff1} depends on the dimension $d$. Indeed, persistence becomes `easier' with increasing dimension. For $d\geq 5$, persistence even becomes the rule: we will show in Section~\ref{example:dgeq5} that if $d\geq 5$ and $\chi$ `truly $d$-dimensional', we always have $(U*U^-)B_0^r<\infty$ (and thus, by Corollary~\ref{cor:suff1}, always persistence of $(\oxi_n)$). This will be proven by means of the following result.
  \begin{corollary}[Sufficient persistence conditions based upon characteristic function]\label{cor:suff2}
   Let $\rho:=\E\chi$, $\rho^-:=\E\chi(-\cdot)$, and $\hat{\rho}(z) := \int\exp(ixz)\E\chi(\d x)$ the characteristic function of the probability distribution $\rho:=\E\chi$. For $\varepsilon>0$ small enough, we have that
   \begin{align}\label{eq:inequalities}
  \sup_{s < 1}\int_{B_0^{\varepsilon} }\frac{1}{|1 - s\hat{\rho}(z)|^2}\d z\leq  \int_{B_0^{\varepsilon}} \frac{1}{|1 - \hat{\rho}(z)|^2}\d z 
\leq \int_{B_0^{\varepsilon}} \frac{1}{|1 - \Re\hat{\rho}(z)|^2}\d z.
   \end{align}
And if any of the integrals is finite for some $\varepsilon>0$ and $\Var\|\chi\|\in(0,\infty)$, then the  critical cluster cascade $(\oxi_n)$ generated by $\mathcal{L}(\chi)$ persists.  
 \end{corollary}
\begin{proof}
As $\Re\hat{\rho}$ is continuous and $\Re\hat{\rho}(0) = 1$, we may pick $\varepsilon>0$ so small that $\Re\hat{\rho}(z)\in\[0,1\]$ for all $z \in B_{0}^\varepsilon$. In this case, we have for $z \in B_{0}^\varepsilon$ and $s\leq 1$
 \begin{equation}\label{eq:proof:inequalities}
 |1-s\hat{\rho}(z)|^2 \geq  |1-\hat{\rho}(z)|^2 \geq (1-\Re\hat{\rho}(z))^2.
 \end{equation}
 The inequalities in \eqref{eq:inequalities} follow from \eqref{eq:proof:inequalities}.

Let $U$ respectively~$U^-$ be the expected occupation measure of a random walk on $\R^d$ with step-size distribution $\rho$ respectively $\rho^-$.
Note that $U*U^- = \sum_{n=0}^\infty\sum_{k = 0}^\infty\rho^{*n}*(\rho^-)^{*k}$. Following exactly the lines of the first part of the proof of Theorem 9.4 (recurrence criterion for random walks) in \cite{kallenberg02}, we find that, for $r >0$,
\begin{equation*}
(U*U^-)B_0^r\leq c' \sup_{s < 1} \int_{B_0^{\sqrt{d}/r}}\frac{1}{1-s\hat{\rho}(z)}\frac{1}{1- s\hat{\rho}_-(z)}\d z
\end{equation*}
 for some finite constant $c'>0$ (depending on $d$ and $r$).
Note that $\hat{\rho}_- = \overline{\hat{\rho}}$ so that the denominator of the integrand becomes
$ (1-s\hat{\rho}(z))(\overline{1-s\hat{\rho}(z)}) = |1-s\hat{\rho}(z)|^2$ and we obtain
\begin{equation*}
(U*U^-)B_0^r\leq c' \sup_{s < 1} \int_{B_0^{\sqrt{d}/r}}\frac{1}{|1-s\hat{\rho}(z)|^2}\d z
\end{equation*}
\change{Thus, taking $r > \varepsilon/\sqrt{\delta}$, we find that $(U*U^-)B_0^r$ is finite whenever one of the integrals in \eqref{eq:inequalities} is finite. Persistence of the critical cluster cascade then follows from Corollary~\ref{cor:suff1}.}
   \end{proof}

\begin{corollary}[Necessary persistence condition]\label{cor:nec}
Let $\Var\|\chi\|\in(0,\infty)$. If the random walk generated by the distribution $\E\chi$ is recurrent, then the critical cluster cascade $(\oxi_n)$ generated by $\mathcal{L}(\chi)$ extinguishes (i.e., it converges weakly to the void point process).
\end{corollary}
\begin{proof}
We show that under the recurrence assumption, $\P\{\oeta_\infty B_0<\infty\} = 0$. Then extinction follows from Theorem~\ref{thm:persistence}.
By Proposition~\ref{prop:direct_construction_of_infinite_palm_tree}, it suffices to show that, under the recurrence assumption, 
 \begin{align}
\ochi^0_{\infty,0}B_0^r  + \sum_{n =0 }^\infty\(\int\ochi^{x}_{\infty,n} B_0^r\beta_{1,n}^{\zeta^-_{n}}(\d x)  +  \beta^{\zeta^-_{n}}_{0,n}B_0^r\) = \infty\quad \text{a.s.}\label{eq:inf}
\end{align}
We observe that the random measure in \eqref{eq:altpalmtree_measure} counts points from the random walk $(\zeta_n^-)$ (the `infinite backward spine') defined in \eqref{eq:rw}. For its step-size distribution we obtain
\begin{align}
\P\left\{\int x \beta_{0,n}(x)\in B\right\}
&\stackrel{\eqref{eq:parentsibling}}{=} \P\left\{\int x \eta_0^{(1)}(x)\in B\right\}\nonumber\\
&\stackrel{\eqref{eq:uetal}}{=} \E\int 1\left\{ \int x \theta_{-y}\chi_0(\d x) \in B\right\}\chi_1(\d y)\nonumber\\
&{=} \E\int 1\left\{ \int x \delta_{-y}(\d x) \in B\right\}\chi_1(\d y)\nonumber\\
&= \E\int 1\{ -y \in B\}\chi_1(\d y)\nonumber\\
&= \E \chi( -B)\nonumber\\
&= \rho^-B,\quad B\in\mathcal{B}(\R^d).\label{eq:rho-}
\end{align} 
  If the random walk generated by $\rho$ is recurrent, then also the random walk generated by $\rho^-$ is recurrent. So \eqref{eq:altpalmtree_measure} a.s.\  observes infinitely many points of the infinite backward spine in $B_0^r$ so that \eqref{eq:inf} holds and therefore $\P\{\oeta_\infty B_0^r <\infty\} = 0$.
\end{proof}

\section{Examples}\label{sec:examples} In this section, we give examples for critical clusters $\chi$ and their critical cluster cascades. \change{For the sake of giving elementary examples, we will also consider some nondiffuse clusters.}
\subsection{Deterministic clusters} Let $\chi:= \delta_{x_0}$ for some $x_0\in\R^d\setminus\{0\}$ so that 
$\ochi_n = \sum_{k = 0}^n \delta_{k{x_0}}.
$ Because ${x_0}\neq 0$,  $\oxi_n B_0^r\leq \lceil2r/|x_0|\rceil \okappa^r_n,\ r>0,\ n\in\N$. Therefore, 
$(\oxi_n B_0^r)_n$ is uniformly integrable (because $( \okappa^r_n)_n$ is uniformly integrable; see Lemma~\ref{lemma:ui}). Thus, $\E \oxi_\infty = c\lambda^d$ and the critical cluster cascade $(\oxi_n)$ persists. We treat the case $x_0 = 0$ in the next example.

\subsection{Clusters without displacements}\label{sec:clusters_without_displacements}
Let $\chi:= Y\delta_0$ for some $\N_0$-valued random variable $Y$ with $\E Y=1$ (possibly $Y\equiv 1$).  Denote by $(Z_n)$ the critical Galton--Watson process generated by $\mathcal{L}(Y)$ (with $Z_0 := 1$).  Then $\ochi_n\eqd \delta_0 \sum_{k=0}^nZ_k$, and 
\begin{align*}
\E\okappa_n^r 
&= \frac{1}{n+1}\int\P\left\{\delta_{0} B_x^r\sum_{k=0}^nZ_k>0\right\}\d x\\
&=\frac{1}{n+1}\int_{B_0^r}\P\left\{\sum_{k=0}^nZ_k > 0\right\}\d x\\
& = \frac{1}{n+1}\int_{B_0^r}\P\{Z_0 > 0\}\d x\\
&= \frac{\lambda^dB_0^r}{n+1}\to0,\quad n\to\infty, r>0.
\end{align*}
Thus, by Theorem~\ref{thm:kappa}, the critical cluster cascade extinguishes.

\subsection{Clusters consisting of exactly one point a.s.\ }\label{subsec:1ptclusters}
Let $\chi:= \delta_X$ for some random variable $X$ on $\R^d$ with distribution $\rho$($=\E\chi$). Note that $\Var \|\chi\| = 0$. We will show that in this case the critical cluster cascade $(\oxi_n)$ persists iff $RW(\rho)$, the random walk generated by $\rho$, is transient:\\ \par
Clearly, $\ochi_n = \sum_l^n\delta_{S_l}$, where $S_0:= 0^d$ and $S_l:= S_{l-1} + X_l,\ l\in\N$, $X_1,X_2,\dots\iidsim \rho$. So $(\ochi_\infty)$ is the occupation measure of $RW(\rho)$. Let $\ochi_n^{\pm}:=\ochi_n + \ochi^-_n - \delta_0$, where $\ochi^-_n$ denotes the occupation measure of the first $n$ steps of a $RW(\rho^-)$ with $\rho^-:=\mathcal{L}(-X)$, independent of  $\ochi_\infty$. Then
\begin{align*}
\E\|\ochi^{r,k}_n\|  &=\E\int1\{\ochi_nB_x^r\leq k\}\ochi_n(\d x)\\
&\geq\E\int1\{\ochi_\infty B_x^r\leq k\}\ochi_n(\d x)\\
&=\sum_{l =0}^{n}\P\{\ochi_\infty B_{S_l}^r\leq k\}\\
&\geq\sum_{l =0}^{n}\P\{\ochi^{\pm}_\infty B_{S_l}^r\leq k\}\\
&=(n+1)\P\{\ochi^{\pm}_\infty B_{S_0}^r\leq k\},\quad r,k>0.
\end{align*}
Consequently,
\begin{align*}
\P\{\oeta_\infty B_0^r<\infty\} 
\stackrel{\mathrm{Lemma\,\ref{lemma:palmtruncation}\emph{\ref{lemma:palmtruncation:limit}}}}{=}
\lim_{k\to\infty}\lim_{n\to\infty}\frac{\E\|\ochi^{r,k}_n\|}{n+1} &\geq  \P\{\ochi_\infty^{\pm} B_{0}^r< \infty\} = \P\{\ochi_\infty B_{0}^r< \infty\} = 1,
\end{align*}
whenever the random walk generated by $\rho$ is transient. Therefore, by Theorem~\ref{thm:persistence}, if the random walk $(S_n)$ generated by $\rho$ is transient, the critical cluster cascade $(\oxi_n)$ persists.

\medskip
On the other hand, note that 
\begin{align*}
\E\|\ochi^{r,k}_n\|  
= \sum_{l = 0}^n\P\{\ochi_nB_{S_l}^r\leq k\}
= \sum_{l = 0}^n\P\{(\ochi_{n-l}  + \ochi^-_{l} -\delta_0)B_{0}^r \leq k\}
&\leq \sum_{l = 0}^n\P\{ \ochi^-_{l}B_{0}^r \leq k\},\quad r,k>0.
\end{align*}
If $\rho$ generates a recurrent random walk, so does $\rho^-$. Consequently $\P\{ \ochi^-_{l}B_{0}^r \leq k\}$ is a zero sequence in $l$ for all $k>0$ and therefore
\begin{equation*}
\lim_{n\to\infty}\frac{1}{n+1}\E\|\ochi^{r,k}_n\| \leq \lim_{n\to\infty}\frac{1}{n+1}\sum_{l = 0}^n\P\{ \ochi^-_{l}B_{0}^r \leq k\} =0
\end{equation*}
as the successive partial averages of a zero sequence converge to zero. So, from Lemma~\ref{lemma:palmtruncation}\ref{lemma:palmtruncation:limit}, we find that $\P\{\oeta_\infty B_0^r < \infty\} = 0$. Thus, by Theorem~\ref{thm:persistence}, we find that the critical cluster cascade $(\oxi_n)$ extinguishes whenever the random walk $(S_n)$ generated by $\rho$ is recurrent.

\subsection{Symmetric $\alpha$-stable cluster intensities}\label{ex:alphastable}
Let $\chi$ a point process on $\R^d$ with intensity $\rho := \E\chi$ a probability measure, following an $\alpha$-stable distribution with characteristic function $\R\ni s\mapsto\exp\{-|s|^\alpha\}$ and $\Var\|\chi\| \in(0,\infty)$. Then we have for all $\varepsilon > 0$ that
\begin{equation}
 \int_{B_0^{\varepsilon}} \frac{1}{|1 - \exp(-|z|^\alpha)|^2}\d z\leq c'\int_0^{\varepsilon} \frac{ s^{d - 1}}{|1 - \exp(- s^\alpha)|^2}\d  s\label{eq1},
\end{equation}
where we change to polar coordinates and use that for all $d$ the modulus of the functional determinant of the transformation is bounded by $s^{d-1}$. Finally, we note that
\begin{equation*}
\frac{ s^{d - 1}}{|1 - \exp(- s^\alpha)|^2} \sim  s^{d - 1 - 2\alpha}, \quad  s\downarrow 0,
\end{equation*}
so that Eq.~\eqref{eq1} is finite when $ d - 1- 2\alpha > -1$.
From Corollary \ref{cor:suff2} if follows that, for $\alpha < d/2$, the corresponding $(\oxi_n)$ persists.  In particular, if $\alpha = 1$ (symmetric Cauchy distribution), $(\oxi_n)$ persists for $d\geq 3$. And if $\alpha =2$ (normal distribution), then $(\oxi_n)$ persists for $d\geq 5$. (We will show in Section~\ref{example:dgeq5} that in fact, for $d\geq5$, the critical cluster cascade $(\oxi_n)$ persists for all `truly $d$-dimensional' critical cluster distributions.)

\subsection{Critical Hawkes processes}\label{example:hawkes}
Hawkes processes, as presented in \cite{hawkes74}, are Poisson cluster point processes on $\R$, where the clusters consist of outgrown subcritical branching random walks generated by finite Poisson processes (with  $\E\|\chi\|\in(0,1)$ and $\E\chi\R_- =0$). In \cite{bremaud01}, a limit construction is considered where the immigration (respectively cluster center) intensity is $\delta c$ for some c>0, and the reproduction mean is $1-\delta$. Letting $\delta\downarrow0$, Theorem 1 in the above paper gives sufficient conditions for local uniform integrability. In the following, we analyze these conditions in our framework.

First of all note that for any distribution $F$ on $\R$
\begin{align*}
| 1- \hat{F}(z)|
&\geq \Re(1 - \hat{F}(z))\\
& = \int 1 - \cos(x z) \d F(x)\\
&\geq \int_{-|z|^{-1}}^{|z|^{-1}} 1 - \cos(x z) \d F(x)\\
&\geq \int_{-|z|^{-1}}^{|z|^{-1}} \frac{(xz)^2}{3} \d F(x)\\
&\geq z^2\int_{-|z|^{-1}}^{|z|^{-1}} \frac{x^2}{3} \d F(x)
\end{align*}
Now consider a probability measure $F$ supported on $\R_+$ s.t. $\overline{F}(x):= 1- F(x) \sim x^{-\alpha}L(x)$ for $\alpha \in(0, 0.5)$ and $L$ slowly varying. For such an $F$, we obtain
\begin{align*}
\int_{-|z|^{-1}}^{|z|^{-1}} {x^2} \d F(x) 
& =\int_{0}^{|z|^{-1}} x^2 \d F(x) \\
&= -\int_{0}^{|z|^{-1}} x^2 \d \overline{F}(x)\\
&= -\[x^2 \overline{F}(x)\]_{x=0}^{|z|^{-1}} + 2\int_{0}^{|z|^{-1}} x \overline{F}( x)\d x\\
&=s |z|^{-2}  \overline{F}(|z|^{-1}) + 2\int_{0}^{|z|^{-1}} L(x)x^{-\alpha+1}\d x\\
&\sim |z|^{-2}  L(|z|^{-1})|z|^{\alpha} + 2\int_{0}^{|z|^{-1}} L(x)x^{-\alpha+1}\d x\\
&=  |z|^{\alpha-2}  L(|z|^{-1}) + 2 (2 - \alpha)^{-1} |z|^{\alpha - 2 } L(|z|^{-1}) \\
& = |z|^{\alpha-2}  L(|z|^{-1}){(4-\alpha)}/{(2-\alpha)}.
\end{align*}
So, as $z\to\infty$, we have
\begin{align*}
 \frac{1}{|1 - \hat{F}(z)|^2} 
 \sim c |z|^{-2\alpha}  L(|z|^{-1})^2.
 \end{align*}
 Consequently, from Corollary~\ref{cor:suff2}, we find that $(\oxi_n)$ persists if $\alpha < 0.5$---thus retrieving the result on the existence of critical Hawkes processes in Theorem 1 of \cite{bremaud01}---without the technical condition on the behavior of $F$ near 0.

\subsection{Poisson clusters}\label{subsec:poisson}
If $\chi$ is a Poisson processs, we can write it as $\chi = \sum_{i = 1}^Y\delta_{X_i}$ with $\{Y, X_1, X_2, \dots\}$ independent, $Y\sim\Pois(1)$, and $X_1, X_2, \dots \iidsim \rho(:=\E\chi)$. The cascade construction is actually very similar to a Hawkes process; in fact, in Example 6.3(c) of \cite{daley03I}, such constructions are in fact called a Hawkes process. One can show that in this case, the parent/siblings process $(\beta_{0,n}, \beta_{1,n})$ from \eqref{eq:parentsibling} is particularly simple, namely
\begin{equation*}
(\beta_{0,n}, \beta_{1,n})\eqd \( \delta_{-X_0}, \theta_{-X_0}\chi\) 
\end{equation*}
with $X_0\sim\rho$, independent of $\chi$. 
The direct construction of the infinite Palm tree in Proposition~\ref{prop:direct_construction_of_infinite_palm_tree} simplifies in an analogous way so that in the Poisson case we simply obtain
\begin{equation}
\oeta_\infty B_0^r \eqd \sum_{n=0 }^\infty \int\ochi^{x}_{\infty,n} B_0^r\zeta^-_{n}(\d x),\label{eq:opi_simple}
\end{equation}
with $(\zeta^-_n)$ the random walk generated by $\rho^- := \E\chi(-\cdot)$; see~\eqref{eq:rho-}. I.e., in the Poisson case, we simply have $\E \oeta_\infty= U*U^-$, and one obtains the sufficient persistence condition from Corollary~\ref{cor:suff1} even more directly.

\subsection{Extinction for dimensions $d=1,2$}\label{example:d=1,2}
Let $\Var\|\chi\|\in(0,\infty)$, $\rho:=\E\chi$ and RW$(\rho)$ the random walk generated by $\rho$.
If $d= 1$ and $\int x\rho(\d x) = 0$, then RW($\rho$) is recurrent; see, e.g., Theorem 9.2 in \cite{kallenberg02}. Thus, by Corollary~\ref{cor:nec}, the critical cluster cascade $(\oxi_n)$ extinguishes. In the case $d=2$, if in addition to the zero mean we have $\int |x|^2\rho(\d x) <\infty$, then, by the same arguments, we also obtain extinction.

\subsection{Persistence for $d\geq 5$} 
\label{example:dgeq5}
If the effective dimension of $\rho:=\E\chi$ (the dimension of the linear subspace spanned by the support of $\rho$) is greater or equal than 5, then the critical cluster cascade $(\oxi_n)$ persists. Indeed, arguing as in the proof of Theorem 9.8 in \cite{kallenberg02} (on transience of random walks with effective dimension $d\geq5$) we find that for all dimensions $\d\in\N$, there are constants $\delta, c>0$, such that
$\big|1-\hat{ \rho}(t)\big|  \geq c|t|^2$ for $ t\in B_0^\delta(d):=B_0^\delta(\subset \R^d).
$ So 
\begin{equation*}
\int_{B_0^\delta(d)}\frac{1}{|1-\hat{\rho}(t)|^2}\d t\leq c\int_{B_0^\delta(d)}|t|^{-4}\d t= c'\int_{0}^\delta r^{-4}r^{d-1}\d r=c'\int_{0}^\delta r^{d-5}\d r,
\end{equation*}
which is finite if $d \geq 5$. So, by Corollary~\ref{cor:suff2}, the critical cluster cascade $(\oxi_n)$ persists.

 \section{Outlook}\label{sec:outlook}
 The main result of the paper are the equivalent formulations of persistence in Theorem~\ref{thm:persistence}. They yields various sufficient and necessary conditions that enable us to present numerous examples. However, the work presented is incomplete in three aspects:

Firstly, we only find the relatively weak necessary condition for persistence in Corollary~\ref{cor:nec}. Therefore, we are only able to give a few examples for extinction; see Section~\ref{example:d=1,2}. For instance, we give no example for extinction for $d=3$ and $d=4$. (Note that we show in Section~\ref{example:dgeq5} that for $d\geq5$ persistence is guaranteed.) And in the Hawkes process context of Section~\ref{example:hawkes} (where $d=1$), we actually know from Proposition 1 in \cite{bremaud01} that if the displacement mean is finite, i.e. $\int x \E\chi(\d x)<\infty$, then persistence is not possible. We were not able to retrieve this necessary persistence condition in our framework.
 
 Secondly, we did not include a systematic comparison between clusters that generate a persistent critical cluster cascade and clusters that are `stable' in the sense of \cite{kallenberg17}. One can show that the first notion implies the second. Furthermore, the notions do not coincide: For the simplest example, consider $\chi:=\delta_0$. This critical cluster is obviously `stable' in the sense of \cite{kallenberg17}.  The corresponding critical cluster cascade, however, is not persistent because in this case $\oxi_nB = n\mu_nB\to 0$ a.s.\  as $n\to\infty$.
 
 \change{Finally, we conjecture that for any possible limit process $\oxi_\infty$ of a  critical cluster cascade generated by a critical cluster distribution $\mathcal{L}(\chi)$, there exists a critical cluster field $(\chi^x)$ such that 
$ \oxi_\infty = \int\chi^x\oxi_\infty(\d x),$ a.s. That is, we think that in the persistent case, the process $\oxi_\infty$ can be represented as a `pathwise solution' of the critical cluster field $(\chi^x)$. Though this seems to be clear intuitively (as the immigrants die out so that all observed particles will have a parent and all branching random walks are outgrown so that the potential offspring of all particles will be included), we were not able to get closer to this conjecture than with the following $L^1$-convergence statement:
 \begin{proposition}\label{prop:cluster_invariance}
Let $(\oxi_n)$ be a critical cluster cascade and $(\chi^x)$ the respective critical cluster field. 
Then, for any bounded Borel set $B\subset\R^d$,
\begin{equation}
\lim_{n\to\infty}\E\left|\oxi_n B -\int\chi^x B\oxi_n(\d x) \right|\to 0.\label{eq:cluster_invariance}
\end{equation}
\end{proposition}
\begin{proof}
Applying the cluster field $(\chi^x)$ to the $n$-th process of the critical cluster cascade corresponds to substituting each point by its (potential) children points. I.e.,
\begin{align*}
\int\chi^x \oxi_{n}(\d x) &= \int (\ochi_{n+1}^x - \chi^x_0) \mu_{n}(\d x) \\
&= \int (\ochi_n^x - \delta_x) \mu_{n}(\d x) + \int \chi^x_{n+1}\mu_{n}(\d x)\\
& = \oxi_n - \mu_n + \int \chi^x_{n+1}\mu_{n}(\d x).
\end{align*}
Consequently, for all $B\in\mathcal{B}_b(\R^d)$,
\begin{align*}
\left|\oxi_n B -\int\chi^x B\oxi_n(\d x) \right|
= \left|\mu_n B -\int\chi^x_{n+1} B\mu_n(\d x)\right|
\leq \mu_n B + \int \chi^x_{n+1} B\mu_n(\d x)
\end{align*}
The right hand side has expected value $2\lambda^dB/(n+1)$ which converges to zero.
\end{proof}
One might call property~\eqref{eq:cluster_invariance} `{pathwise cluster invariance}': for $n$ large enough, applying `clustering' to all points of $\xi_n$ does not change the process. Note however, that though the clusters $\chi^x$ are independent over $x\in\R^d$, the points of $\xi_n$ and of the clusters are not---in contrast to `normal' clustering, where clusters are independent of the realization of the argument process. Perhaps a better rewording of \eqref{eq:cluster_invariance} might be that $\xi_n$ `solves the cluster field $(\chi_x)$' for large $n$.
}
\section{Acknowledgements}
The detailed critical comments of two anonymous referees together with remarks of Paul Embrechts on a first version of the paper improved overall presentation considerably. Furthermore, we thank Eugen Jost for supplying the illustrations.

\appendix
\section{Proofs}\label{appendix:proofs}
\begin{proof}[Proof of Lemma \ref{lemma:monotonicity}]\label{proof:lemma:monotonicity}
Finiteness of $\E\okappa_n^r$ immediately follows from local finiteness of the intensity of $\oxi_n$:
\begin{align*}
\E\okappa_n^r &=\E \int1\{\overline{\chi}^x_nB_0^r >0\}\mu_n(\d x) \leq\E \int\overline{\chi}^x_nB_0^r\mu_n(\d x) =  \E\overline{\xi}_nB_0^r= c\lambda^dB_0^r<\infty, \quad n\in\N_0.\\
\end{align*}For monotonicity, we first note that for all $x\in\R^d$
\begin{align}
&\frac{1}{n+1}\P\{\ochi_nB_x^r>0\} -\frac{1}{n}\P\{\ochi_{n-1}B_x^r>0\} \nonumber\\
&=\frac{1}{n+1}\P\{\ochi_nB_x^r>0\} -\(\frac{1}{n+1} +\frac{1}{n(n+1)}\)\P\{\ochi_{n-1}B_x^r>0\}\nonumber\\
& =\frac{1}{n+1}\P\{\chi_nB_x^r>0, \ochi_{n-1} B_x^r = 0\} - \frac{1}{n(n+1)}\P\{\ochi_{n-1}B_x^r>0\}\nonumber\\
& = \frac{1}{n(n+1)}\(n\P\{\chi_nB_x^r>0, \ochi_{n-1} B_x^r = 0\} - \sum_{k = 0}^{n-1}\P\{\chi_{k}B_x^r>0, \ochi_{k-1}B_x^r =0\}\), \label{eq:diff}
 \end{align}
 where we set $\chi_{-1}\R^d := 0$. For the term in the bracket of \eqref{eq:diff}, we obtain
  \begin{equation}
 \sum_{k = 0}^{n-1}\P\{\chi_nB_x^r>0, \ochi_{n-1} B_x^r = 0\}  - \P\{\chi_{k}B_x^r>0, \ochi_{k-1}B_x^r =0\}.\label{eq:sum}
 \end{equation}
 To prove that \eqref{eq:sum}---and therefore \eqref{eq:diff}---is nonpositive, it suffices to show that all summands in \eqref{eq:sum} are nonpositive. Observing that both summands are values of the sequence $(\P\{\chi_{k}B_x^r>0, \ochi_{k-1}B_x^r =0\})_{k\in\N_0}$, it furthermore suffices to show that this sequence is nonincreasing. Indeed, we have
 \begin{align}
 &\int\P\{\chi_{k+1}B_x^r>0, \ochi_{k}B_x^r =0\}\d x\nonumber\\
 & =  \int\E1\left\{\int \chi_k^uB_x^r\chi(\d u)>0,\ \delta_0B_x^r +\int \ochi_{k-1}^uB_x^r\chi(\d u)=0 \right\}\d x\nonumber\\
 &\leq\int\E1\left\{\int \chi_k^uB_x^r\chi(\d u)>0,\ \int \ochi_{k-1}^uB_x^r\chi(\d u)=0 \right\}\d x.\label{eq:indicator}
 \end{align}
 Note that,  for fixed $x\in\R^d$, the event in the indicator function in \eqref{eq:indicator} can be rewritten in the following way:
 \begin{align*}
 &\left\{\int \chi_k^uB_x^r\chi(\d u)>0,\ \int\ochi_{k-1}^uB_x^r\chi(\d u)=0\right\}\\
 &= \left\{1\{\chi_k^uB_x^r>0\} > 0 \text{ for some }  u \text{ with }\chi(\d u) = 1,\ 1\{\ochi_{k-1}^uB_x^r=0\} > 0 \text{ for all } u\text{ with }\chi(\d u) = 1\right\}\\
 &=\left\{\int 1\{\chi_k^uB_x^r>0, \ochi_{k-1}^uB_x^r=0\} \chi(\d u) > 0\right\}.
\end{align*}
Consequently, we get the following upper bound for the indicator function in \eqref{eq:indicator}:
\begin{align*}
&1\left\{\int \chi_k^uB_x^r\chi(\d u)>0,\ \int \ochi_{k-1}^uB_x^r\chi(\d u)=0 \right\}\\
&=1\left\{\int 1\{\chi_k^uB_x^r>0,\  \ochi_{k-1}^uB_x^r=0 \}\chi(\d u)>0\right\}\\
&\leq \int 1\left\{\chi_k^uB_x^r>0,\  \ochi_{k-1}^uB_x^r=0 \right\}\chi(\d u).
\end{align*}
  Plugging this upper bound into \eqref{eq:indicator}, we obtain
   \begin{align*}
  \int\P\{\chi_{k+1}B_x^r>0, \ochi_{k}B_x^r =0\}\d x
&\leq 
  \int\E\int 1\{\chi_k^uB_x^r>0,\  \ochi_{k-1}^uB_x^r=0 \}\chi(\d u)\d x\\
  &=\int\int\P\{ \chi_k^uB_x^r>0,\ \ochi_{k-1}^uB_x^r=0 \}\E\chi(\d u)\d x\\
  &=\int\int\P\{ \chi_k^uB_x^r>0,\ \ochi_{k-1}^uB_x^r=0 \}\d x\E\chi(\d u)\\
      &= \int\int 1\{\chi_{k}B_{x-u}^r>0, \ochi_{k-1}B_{x-u}^r =0\}\d xE\chi(\d u)\\
    &= \int\int\P\{\chi_{k}B_{x}^r>0, \ochi_{k-1}B_{x}^r =0\}\d x\E\chi(\d u)\\
  &= \int\P\{\chi_{k}B_x^r>0, \ochi_{k-1}B_x^r =0\}\d x {\E\|\chi\|}.
 \end{align*}
 So, the sequence $(\int\P\{\chi_{k}B_x^r>0, \ochi_{k-1}B_x^r =0\}\d x)_{k}$ is nonincreasing. Therefore, the summands in \eqref{eq:sum} are nonpositive and, consequently, also \eqref{eq:diff} is nonpositive. So $(\E\okappa_n^r)_{n}$ is nonincreasing.\\
 For the last statement in the lemma, choose $r_0>0$ such that $\lim_{n\to\infty} \E\okappa_n^{r_0} >0$. The limit is also positive for $r\geq r_0$ because $\E\okappa_n^r$ is nondecreasing in $r$. If $r\leq r_0$, cover the ball $B^{r_0}_0$ with finitely many, say $m$, balls of the form $B_{x_k}^r$. We obtain
\begin{align*}
\E\okappa_n^{r_0}
= \frac{c}{n+1}\int\P\{\overline{\chi}_nB_x^{r_0} >0\}\d x &\leq \frac{c}{n+1}\int\P\{\overline{\chi}_n\cup_{k=1}^mB_{x+x_k}^r >0\}\d x\\
& \leq \sum_{k = 1}^m\frac{c}{n+1}\int\P\{\overline{\chi}_nB_{x}^r >0\}\d x\\
& = m\E\okappa_n^r.
\end{align*}
So, if $\E\okappa_n^{r_0}$ has a strictly positive limit, also $\E\okappa_n^r$, $0<r\leq r_0$, will have a strictly positive limit.
\end{proof}

\begin{proof}[Proof of Proposition~\ref{prop:direct_construction_of_infinite_palm_tree}]\label{proof:prop:direct_construction_of_infinite_palm_tree}
Recall from \eqref{eq:L} that $L_{n+1}-L_{n}\in\{0,1\}$ for all $n$. If $L_{n+1}-L_{n} = 0$, the construction given in \eqref{eq:fwbw_recursion} makes a forward step; if $L_{n+1}-L_{n} = 1$, the construction makes a backward step. First, we will show that there are a.s\ infinitely many forward and backward steps:

Since, $(L_n)$ is nondecreasing, and $L_n\sim\mathrm{Unif}\{0,1,2,\dots n\},\, n\in\N_0,$
\begin{align*}
\P\{\lim_{n\to\infty} L_n <\infty\}
= \lim_{l\to\infty}\lim_{n\to\infty}\P\{L_n \leq l\} =  \lim_{l\to\infty}\lim_{n\to\infty}\frac{l+1}{n+1}=0.
\end{align*}
I.e., $\lim_{n\to\infty} L_n= \infty$ a.s\ and consequently $\{L_n - L_{n-1} = 1 \text{ i.o.}\}$ a.s. Similarly, one can show that $\lim_{n\to\infty} (n - L_n) = \infty$ a.s.\ so that also $\{L_n - L_{n-1}= 0  \text{ i.o.}\}$.

So in the construction of the nondecreasing sequence $\oeta_n^{(L_{n})}$ as defined in \eqref{eq:fwbw_recursion}, we have infinitely many forward steps (first case in \eqref{eq:fwbw_recursion}). I.e., from any possible point, we grow an infinite (forward) cumulative critical branching random walk. Secondly, there are also infinitely many backward steps (second case in \eqref{eq:fwbw_recursion}). I.e., we attach a parent (i.e., a new root) and siblings to each previous root.

Consequently, if we are only interested in the limit $\oeta_{\infty} = \lim_{n\to\infty}\oeta_n^{(L_{n})}$, we may ignore the sequence $(L_n)$, i.e., the decisions between backward and forward steps. Instead, we start with a single point in zero and immediately attach the infinite backward spine of parents and siblings and then attach outgrown branching random walks $\ochi_{\infty,n}^x$ to the zero point and all sibling points which gives representation \eqref{eq:altpalmtree}.
\end{proof}

\begin{proof}[Proof of Lemma~\ref{lemma:bound}]\label{proof:lemma:bound}
We aim to show that $\E \oeta_\infty B_0^r <\infty$ whenever $(U*U^-) B_0^{2r} <\infty$. In the case with Poisson clusters this follows immediately because, in this case, straightforward calculations show that $\E \oeta_\infty  = U*U^- $; see Section~\ref{subsec:poisson}.

The general case demands more argumentation: consider the representation of $\oeta_\infty$ in Proposition~\ref{prop:direct_construction_of_infinite_palm_tree}. Taking the expectation of the $n$-th summand in \eqref{eq:altpalmtree_measure} we obtain
\begin{align*}
\E \int \ochi_{\infty,n}^{x}\beta^{\zeta^-_{n}}_{1,n}(\d x) + \E \beta_{0,n}^{\zeta^-_{n}}
&=\E\int\int \ochi^x_{\infty,n}  \beta_{1,n}^{z}(\d x)\delta_{\zeta^-_{n}}(\d z)  + \E \delta_{\zeta^-_{n+1}}\\
&=\int\int \E\ochi_{\infty,n}(\cdot - x)^r \E\beta^z_{1,n}(\d x)\E\delta_{\zeta^-_{n}}(\d z)  + (\rho^-)^{(n+1)*},
\end{align*}
where in the first summand, we used that $\zeta^-_{n}$ and all the $\ochi^x_\infty$ are independent of each other as well as of $(\beta^x_{0,n},\beta^x_{1,n})$, and in the second term, we denote $\rho^-:=\E\chi(-\cdot)$, the step-size distribution of the random walk $(\zeta^-_n)$; see \eqref{eq:rho-}. Noting that $\E\ochi_\infty = \sum_{k=0}^\infty\rho^{k*}=:U$, with $\rho:=\E\chi$ and $\rho^{0*}:= \delta_0$, we finally find for the expectation of the $n$-th summand in \eqref{eq:altpalmtree_measure}:  
\begin{align*}
\int \int U(\cdot-x)\E\beta^z_{1,n}(\d x)(\rho^-)^{n*}(\d z) + (\rho^-)^{(n+1)*}
&=\int  (U*\E\beta^z_{1,n}) (\rho^-)^{n*}(\d z) + (\rho^-)^{(n+1)*}\\
&=\int  (U*\E\beta_{1})(\cdot-z) (\rho^-)^{n*}(\d z) + (\rho^-)^{(n+1)*}\\
&= (U *\E\beta_{1}* (\rho^-)^{n*}) + (\rho^-)^{(n+1)*}.
\end{align*}
So we derive expected value of \eqref{eq:altpalmtree_measure} by summing the latter formula over $n\geq 0$ and also taking the first term of \eqref{eq:altpalmtree_measure} into account:
\begin{align}
U + (U *\E \beta_1 *  U^-)+ (U^--\delta_0)
&= (U + U^- -\delta_0) + (\E \beta_1 * U * U^-).\label{eq:renewal}
\end{align}
We show that if $\Var\|\chi\| \in(0,\infty)$, then \eqref{eq:renewal} is locally finite whenever
$ U * U$ is locally finite:
First of all, if $U * U^-$ is locally finite, then $U$ as well as $U^-$ are locally finite. So \eqref{eq:renewal} is locally finite iff $(\E \beta_1 * U * U^-)$ is locally finite. Furthermore, note that
\begin{align*}
\E\|\beta_1\|
= \E\|\eta^{(1)}_1-\delta_0\|
= \E\int \|\chi\|\chi(\d x)  -1
= \Var\|\chi\|,
\end{align*}
 so that $F := \E \beta_1/{\Var \|\chi\|}$ defines a probability measure on $\R^d$. Clearly, the measure $(\E \beta_1 * U * U^-)$ is locally finite iff $
 (F * U * U^-)$ is locally finite. The measure $(F * U * U^-)$ has a simple interpretation:
 
\change{On each particle of a random walk $(\zeta^-_n)$ generated by $\rho^-$, we attach (independent) random walks $(\zeta_{n,k})_k$  generated by $\rho$. Then $U * U^-$ denotes the expected occupation measure of this object. Consequently, 
$F*U * U^-$ can be interpreted as having the same construction, where the first point $\zeta^-_0$ of $(\zeta^-_n)$ is `delayed' by the distribution $F$. That is, all the points are shifted by a random variable following distribution $F$.} Using the Markov property of the component processes $(\zeta^-_n)$ and $(\zeta_{n,k})_k$, one can show that for all $r >0$
$(U * U^-)B_x^r\leq (U * U^- )B_0^{2r},\ x\in\R^d.
$ Consequently,
\begin{equation*}
(F*U * U^-)B_0^r = \int (U * U^-)B_{x}^r F(\d x)\leq\int (U * U^-)B_{0}^{2r} F(\d x) = (U * U^-)B_0^{2r}.
\end{equation*}
\change{
Summarizing, we have shown that
$$
\E\oeta_\infty B_0^r= (U + U^- -\delta_0)B_0^r + (\E \beta_1 * U * U^-)B_0^r\leq c_r + \Var\chi(F*U * U^-)B_0^r
\leq c_r + (U * U^-)B_0^{2r}
$$
with $c_r>0$ some finite constant depending on $r>0$.
Thus, if $(U * U^-) B_0^{2r}$ is finite, then $\E\oeta_\infty B_0^r$ is finite.}
\end{proof}
\newpage
 \section{Figures}\label{appendix:figures}

\vfill
\clearpage
%
\bibliographystyle{apalike}
\bibliography{../../../Bibliographies/diss.biblio.bib}

\end{document}